\documentclass[12pt,psfig,reqno]{amsart}
\usepackage{mathrsfs}
\usepackage[colorlinks,
            linkcolor=cyan,
            anchorcolor=cyan,
            citecolor=cyan
            ]{hyperref}
\usepackage{txfonts}
\usepackage{amscd}
\usepackage{cite}
\usepackage{epsfig}
\usepackage{verbatim}
\usepackage{mathdots}
\usepackage{amssymb}
\usepackage{amsfonts}
\usepackage{amsbsy}
 \usepackage{graphicx}
 \usepackage{epstopdf}
 \usepackage{latexsym,amsmath,amssymb,amsthm,amsfonts,arydshln,enumerate,graphicx}
 \usepackage[usenames]{color} 
 \usepackage{cancel}
 \numberwithin{equation}{section}
 
\newtheorem{prop}{Proposition}[section]
\newtheorem{lem}[prop]{Lemma}

\newtheorem{defi}{Definition}[section]
\newtheorem{cor}[prop]{Corollary}
\newtheorem{thm}[prop]{Theorem}

\newtheorem{exam}[prop]{Example}
\newtheorem{rem}[prop]{Remark}

\begin{document}
\baselineskip=17pt
\title{Spectrality of  a class of infinite convolutions on $\mathbb{R}$}
\author{Sha Wu}
\address{Sha Wu, School of Mathematics, Hunan University, Changsha, 410082, P.R. China}
\email{shaw0821@163.com}
\author{Yingqing Xiao$^*$}
\address{Yingqing Xiao, School of Mathematics and Hunan Province Key Lab of Intelligent Information Processing and Applied Mathematics, Hunan University, Changsha, 410082, P.R. China}
\email{ouxyq@hnu.edu.cn}

\date{\today}
\thanks{The research is supported in part by National Natural Science Foundation of China (Grant Nos. 12071118).
\\
$^*$Corresponding author.}
\begin{abstract}
 Given an integer $m\geq1$. Let $\Sigma^{(m)}=\{1,2, \cdots, m\}^{\mathbb{N}}$ be a symbolic space, and let $\{(b_{k},D_{k})\}_{k=1}^{m}:=\{(b_{k}, \{0,1,\cdots, p_{k}-1\}t_{k}) \}_{k=1}^{m}$
be a finite sequence pairs, where integers $| b_{k}| $, $p_{k}\geq2$, $|t_{k}|\geq 1$ and  $ p_{k},t_{1},t_{2}, \cdots, t_{m}$ are pairwise coprime  integers for all $1\leq k\leq m$.  In this paper, we show that
for any infinite word $\sigma=\left(\sigma_{n}\right)_{n=1}^{\infty}\in\Sigma^{(m)}$, the infinite convolution
$$
\mu_{\sigma}=\delta_{b_{\sigma_{1}}^{-1} D_{\sigma_{1}}} * \delta_{\left(b_{\sigma_{1}} b_{\sigma_{2}}\right)^{-1} D_{\sigma_{2}}} *  \delta_{\left(b_{\sigma_{1}} b_{\sigma_{2}}  b_{\sigma_{3}}\right)^{-1}D_{\sigma_{3}}} * \cdots
$$
is a spectral measure if and only if   $p_{\sigma_n}\mid b_{\sigma_n}$  for  all $n\geq2$ and $\sigma\notin \bigcup_{l=1}^\infty\prod_{l}$, where $\prod_{l}=\{i_{1}i_{2}\cdots i_{l}j^{\infty}\in\Sigma^{(m)}: i_{l}\neq j, |b_{j}|=p_{j}, |t_{j}|\neq1\}$.
\end{abstract}
\maketitle
\section{\bf Introduction\label{sect.1}}
\medskip
Let $\mu$ be a Borel probability measure on $\mathbb{R}^{d}$ with compact support $K$. If there exists a countable discrete set $\Lambda\subset \mathbb{R}^{d}$
such that $E(\Lambda)=\{e^{2\pi i<\lambda,x>}:\lambda\in \Lambda\}$ forms an orthogonal basis for the Hilbert space $L^{2}(\mu)$, then $\mu$ is called a \emph{spectral measure} and $\Lambda$ is called a \emph{spectrum} of $\mu$. We also say  that $(\mu ,\Lambda)$ forms a \emph{spectral pair}. In particular, if  $K$ has positive Lebesgue measure and $\mu$ is the Lebesgue measure on $K$,  then $K$ is called a \emph{spectral set}. The classical example of spectral sets is the unit cube $K=[0, 1]^{d}$, for which the set $\Lambda= \mathbb{Z}^{d}$ serves as a spectrum. Moreover, Fuglede \cite{Fug74} proved that triangles and  disks in the plane are not spectral sets.

The research of spectral measures was originated from Fuglede's  spectral set conjecture \cite{Fug74}. Fuglede conjectured that $K\subset \mathbb{R}^d $ is a  spectral set  if and only if $K$ is a translational tile. We say that  $K$ is a \emph{translation tile} if  there exists a discrete set $J\subset \mathbb{R}^d $ such that
 $$\mathbb{R}^n=\bigcup_{j \in J}(K+j)  \quad \text { and } \quad \mathcal{L}((K+j_1) \cap (K+j_2))=0 ~ \text { for all }j_1 \neq j_2\in J,$$ where $\mathcal{L}(\cdot)$ denotes the Lebesgue measure. There is a large literature concerning the translational tiles (see \cite{LW96,HKL11,LW96-1,LW96-2} and the references therein).
Fuglede's conjecture has attracted many researchers in  the past few decades. It has been proved to be false for $d\geq3$ \cite{KM061,Ta04}, and it is still open in one and two dimensions.  Nevertheless, the conjecture is true if some additional conditions are added to $K$. For examples, {\L}aba  \cite{L00} proved that the conjecture is true for the union of  two disjoint intervals of non-zero length in one dimension. More recently, a significant progress on the conjecture has be achieved by  Lev and Matolcsi \cite{LM19}, they  proved that  the conjecture holds in any dimension for a convex body.  For more discussion on  the conjecture, the reader can refer to \cite{L02, LS01, GL16} etc.

After the original work of Fuglede, the study of spectral measures is also blooming in the field of fractal geometry. The first fractal spectral measure, 1/4-Cantor measure, was discovered by Jorgensen and Pedersen \cite{JP98}.
This opens a new field for the study of orthogonal harmonic analysis of fractal measures including self-similar/self-affine measures and generally Moran measures, see \cite{CL19,CLS20,DHL2014,DHL2013,DHL2019,LD2020,LDL17,FHL15,FY97,LMW-2022-2, De14,LL17,L2014,FHW18,LMW-2022-1} and the references therein.

The present paper is devoted to investigating the spectrality of a class Moran measures on $\mathbb{R}$.
 Let  $\left\{b_{k}\right\}_{k=1}^{\infty}$ be a sequence of  integers with  $| b_{k}|\geq2$ and $\left\{D_{k}\right\}_{k=1}^{\infty}$ be a sequence of  digit sets with $D_{k}\subset \mathbb{Z}$. For a finite subset $E\subset\mathbb{R}$, we define $\delta_E=\frac{1}{\#E}\sum_{e\in E}\delta_e$, where $\#E$ denote the cardinality of $E$ and $\delta_e$ is the Dirac point mass measure at $e$. Write
\begin{equation}\label{eq(2.81)}
\mu_{k}=\delta_{b_{1}^{-1}D_1}\ast\delta_{b_{1}^{-1}b_{2}^{-1}D_2}
\ast\delta_{b_{1}^{-1}b_{2}^{-1} b_{3}^{-1}D_3}\ast\cdots\ast\delta_{b_{1}^{-1}b_{2}^{-1}\cdots b_{k}^{-1}D_k},
\end{equation}
where $ \ast$ is the convolution sign.
We say that $\mu_{k}$ converges weakly to $\mu $ if
$$\lim\limits_{k\rightarrow\infty}\int fd\mu_{k}=\int fd\mu$$
 for all $f\in C_{b}(\mathbb{R})$, where  $C_{b}(\mathbb{R})$ denotes the set of all bounded continuous functions on  $\mathbb{R}$.  It is known\cite{SS00} that if
 \begin{equation}\label{eq(2.82)}
\sum_{k=1}^{\infty} \frac{\max \left\{| d |: d \in D_{k}\right\}}{\left|b_{1} b_{2} \cdots b_{k}\right|}<\infty,
 \end{equation}
 then $\mu_{k}$ converges weakly to $\mu$ with compact support.
The weak limit measure $\mu$  may be written as
\begin{equation}\label{eq(2.8)}\begin{aligned}
\mu :&=\mu_{\{b_k\},\{D_k\}}\\&=\delta_{b_{1}^{-1}D_1}\ast\delta_{b_{1}^{-1}b_{2}^{-1}D_2}
\ast\delta_{b_{1}^{-1}b_{2}^{-1} b_{3}^{-1}D_3}\ast\cdots\ast\delta_{b_{1}^{-1}b_{2}^{-1}\cdots b_{k}^{-1}D_k}\ast\delta_{b_{1}^{-1}b_{2}^{-1}\cdots b_{k+1}^{-1}D_{k+1}} \ast \cdots
\\&=\mu_{k}\ast\nu_{>k}(b_{1}b_{2}\cdots b_{k}\cdot),
\end{aligned}\end{equation}
where
\begin{equation}\label {eq(2.90)}
\nu_{>k}:=\delta_{b_{k+1}^{-1}D_{k+1}}\ast\delta_{b_{k+1}^{-1}b_{k+2}^{-1}D_{k+2}}
\ast\delta_{b_{k+1}^{-1}b_{k+2}^{-1}b_{k+3}^{-1}D_{k+3}}\ast\cdots,
\end{equation}
and the measure $\mu$ is called Moran measure (or called Riesz product measure). Moreover, the Moran measure $\mu$  is supported on a compact set
\begin{equation}\label{eq(0.0)}
K:=K(b_{k}, D_{k}) =\left\{\sum_{k=1}^{\infty} \frac{ d_{k} }{b_{1} b_{2} \cdots b_{k}} : d_{k}\in D_{k},~k\geq1\right\},
 \end{equation}
where the set $K$ is usually called a Moran set.

 The spectrality of Moran measure were first studied by Strichartz\cite{SS00} in 2000.  After that,  many  spectral Moran measures were found in \cite{AFL19, AH14,HH17,Tang2018,WDL2018} etc. For examples, in 2014, An and He \cite{AH14}  studied the  spectrality  of Moran measure $\mu_{\{b_k\},\{D_k\}}$  generated by  an integer sequence $\{b_{k}\}_{k=1}^{\infty}$ with $b_k\geq 2$ and a sequence  of consecutive digit sets $\{D_{k} \}_{k=1}^{\infty}$, where $D_{k}=\{0,1,\cdots,p_{k}-1\}$ with $p_{k}\geq2$. They showed that the Moran measure $\mu_{\{b_k\},\{D_k\}}$ is a spectral measure if  $p_{k}\mid b_{k} $  for all $k\geq1$.
 In 2021, under the condition that $\left\{p_{k} \right\}_{k=1}^{\infty}$  is  a bounded sequecne, the sufficient and necessary conditions for $\mu_{\{b_k\},\{D_k\}}$  to become a spectral measure have been given by Deng and Li\cite{DL21}.

\begin{thm}\label{thm(1)}\cite{DL21}
Let $\mu_{\{b_k\},\{D_k\}}$ be the Moran measure generated by an integer sequence  $\{b_k\}_{k=1}^{\infty}$ with $|b_k|\geq 2$ and a sequence of consecutive digit sets $\{D_k\}_{k=1}^{\infty}$, where $D_{k}=\left\{0,1, \cdots, p_{k}-1\right\}$ with $p_k\geq2$. Assume that the sequence $\left\{p_{k} \right\}_{k=1}^{\infty}$ is bounded. Then $\mu_{\{b_k\},\{D_k\}}$ is a spectral measure if and only if  $p_{k}\mid b_{k} $  for all $k\geq2$.
\end{thm}
 In \cite{LD2020}, Liu and Dong considered the spectrality of a class of more general Moran measure $\mu_{\rho^{-1},\{D_{k}\}}$, where $0<\rho<1$ is real  and $D_{k}=\{0,1, 2,\cdots, N^s-1\}t_{k} $  with a bounded integer sequence $\{t_{k}\}_{k=1}^{\infty}$ and a prime $N\geq2$. They proved that $L^2(\mu_{\rho^{-1},\{D_{k}\}})$ admits an infinite set of orthonormal exponential functions if and only if there are  integers $p, q, r \geq1$ such that $\rho^{-1}=(\frac{Np}{q})^{\frac{1}{r}}$ with $\gcd(Np,q)=1$. Moreover, they demonstrated that the necessary condition for  $\mu_{\rho^{-1},\{D_{k}\}}$  to become a spectral measure is $\rho^{-1}=Np$.

Inspired by the above results, we  consider the spectrality of  the Moran measure $\mu_{\{b_k\},\{D_k\}}$, which is generated by an integer sequence $\{b_k\}_{k=1}^{\infty}$ with $|b_k|\geq 2$  and an integer  sequence of digit sets $\{D_k\}_{k=1}^{\infty}$, where $D_{k}=\left\{0,1, \cdots, p_{k}-1\right\}t_{k}$ with $p_k\geq2$, $| t_{k}|\geq1$.
In order to achieve a convenient statement of our results, we first introduce some notations. Given an integer $m\geq1$. Let $\Sigma^{(m)}=\{1,2, \cdots, m\}^{\mathbb{N}}$ be a symbolic space and
\begin{equation}\label{eq(1.11)}\{(b_{k},D_{k})\}_{k=1}^{m}:=\{(b_{k}, \{0,1,\cdots, p_{k}-1\}t_{k}) \}_{k=1}^{m}
\end{equation}
be a finite sequence pairs, where integers $| b_{k}| $, $p_{k}\geq2$,$|t_{k}|\geq 1$ and  $ p_{k},t_{1},t_{2}, \cdots, t_{m}$ are pairwise coprime  integers for all $k\in\{1,2, \cdots, m\}$. For an infinite word $\sigma=\left(\sigma_{n}\right)_{n=1}^{\infty} \in\Sigma^{(m)}$, we define the infinite convolution
\begin{equation}\label{eq(2.12)}
\mu_{\sigma}=\delta_{b_{\sigma_{1}}^{-1} D_{\sigma_{1}}} * \delta_{\left(b_{\sigma_{1}} b_{\sigma_{2}}\right)^{-1} D_{\sigma_{2}}} *  \delta_{\left(b_{\sigma_{1}} b_{\sigma_{2}}  b_{\sigma_{3}}\right)^{-1}D_{\sigma_{3}}} * \cdots.
\end{equation}
The main purpose of this paper is to give the necessary and sufficient condition for $\mu_{\sigma}$ to be a spectral measures, which can be stated as follows.
 \begin{thm}\label{thm(1.13)}
  Given  an integer $m\geq1$. Let $\Sigma^{(m)}=\{1,2, \cdots, m\}^{\mathbb{N}}$ be the symbolic space and $\{(b_{k},D_{k})\}_{k=1}^{m}$ be defined by \eqref{eq(1.11)}. For $\sigma \in\Sigma^{(m)}$, then the infinite  convolution $\mu_{\sigma}$ defined by \eqref{eq(2.12)} is a spectral measure if and only if   $p_{\sigma_n}\mid b_{\sigma_n}$ for  all $n\geq2$ and $\sigma\notin \bigcup_{l=1}^\infty\prod_{l}$, where $\prod_{l}=\{i_{1}i_{2}\cdots i_{l}j^{\infty}\in\Sigma^{(m)}: i_{l}\neq j, |b_{j}|=p_{j}, |t_{j}|\neq1\}$.
\end{thm}
We remark that Theorem \ref{thm(1.13)} extends the following result  of  Li et al. \cite{LMW-2022-2}  into a more general form.
\begin{thm}\label{thm(1.131)} \cite{LMW-2022-2}
Given an integer $k\geq 1$ and two coprime integers $p, t \geq 2$. Let
$$\Sigma^{(2)}=\{1,2\}^{\mathbb{N}}, \ \ b_{1}=b_{2}=kp, \ \ D_1=\left\{0,1, \cdots, p-1\right\}, \ \ D_2=\left\{0,1, \cdots, p-1\right\}t.
$$
\begin{enumerate}[(1)]
 \item If $k=1$, then for $\sigma\in \Sigma^{(2)}$, the infinite  convolution $\mu_{\sigma}$ is a spectral measure if and only if  $\sigma=2^{\infty}$ or the symbol $``1"$ occurs infinitely many times in $\sigma$.
 \item  If $k\geq 2$, then for all $\sigma\in\Sigma^{(2)}$, the infinite convolution $\mu_{\sigma}$ is a spectral measure.
 \end{enumerate}
\end{thm}
In Theorem \ref{thm(1.131)}, we have  $\sigma\notin \bigcup_{l=1}^\infty\left\{i_{1}i_{2}\cdots i_{l }2^{\infty}\in\Sigma^{(2)}: b_{2}=p, i_{l}=1\right\}$ if  $\mu_{\sigma}$ is a spectral measure.

The necessary condition $p_{\sigma_n}\mid b_{\sigma_n}(n\geq2)$ of Theorem \ref{thm(1.13)} can be
obtained directly  by the following Theorem \ref{thm(1.1)}, which eased some restrictions on $  b_{k},p_{k}$ and $t_{k}$ of  Theorem \ref{thm(1.13)}.
\begin{thm}\label{thm(1.1)}
 Given integer sequences $\{b_k\}_{k=1}^{\infty}$, $\{p_k\}_{k=1}^{\infty}$ and $\{t_k\}_{k=1}^{\infty}$ with $p_{k} , |b_{k}| \geq 2$ and $|t_{k}| \geq 1$. Let $\mu_{\{b_k\},\{D_k\}}$ and $\nu_{>k}$ be defined by \eqref{eq(2.8)}  and \eqref{eq(2.90)},  respectively, where $D_{k}=\left\{0,1, \cdots, p_{k}-1\right\}t_{k}$ and $\left\{p_{k}t_{k}\right\}_{k=1}^{\infty}$ is a bounded sequecne. Suppose that $p_{k}\nmid t_{k+1}$ for all $k\geq1$  and $\mu_{\{b_k\},\{D_k\}}$ is a spectral measure, then $p_{k+1}\mid(b_{k+1}t_{k })$ and $\nu_{>k }$ is a spectral measure for all $k\geq1$.
\end{thm}
\begin{rem}\label{rem4.1}
To get $p_{k+1}\mid(b_{k+1}t_{k })$ in Theorem \ref{thm(1.1)}, the condition $p_{k}\nmid t_{k+1}$ is necessary. For example, let $D_{2k-1}=\{0,1\}$, $D_{2k}=\{0,1,2\}4$ and $b_{1}=12$, $b_{2k+1}=6$, $b_{2k}=2$ for $k\geq 1$. As $p_{2k-1}=2$ and $t_{2k}=4$, one can easily know that $p_{2k-1}\mid t_{2k}$ for $k\geq 1$.  We can rewrite $\mu_{\{b_k\},\{D_k\}}$ as follows form
$$
\begin{aligned}
\mu_{\{b_k\},\{D_k\}}&=\delta_{12^{-1}\{0,1\}}\ast\delta_{12^{-1}2^{-1}\{0,1,2\}4}
\ast\delta_{12^{-1}2^{-1}6^{-1}\{0,1\}}\ast\delta_{12^{-1}2^{-2}6^{-1}\{0,1,2\}4}\ast\cdots \\
&=\delta_{12^{-1}\{0,1,2,3,4,5\}}\ast\delta_{ 12^{-2}\{0,1,2,3,4,5\}}
\ast\delta_{ 12^{-3}\{0,1,2,3,4,5\}}\ast\cdots\\
&:=\mu_{12,\{0,1,2,3,4,5\}}.
\end{aligned}
$$
Together with Theorem \ref{thm(1)}, it shows that $\mu_{\{b_k\},\{D_k\}}$ is a spectral measure.
But it is easy to see that $p_{2k}\nmid(b_{2k}t_{2k-1}) $ since $p_{2k}=3$, $b_{2k}=2$ and  $t_{2k-1}=1$ for $k\geq 1$.
\end{rem}
 The following Theorem \ref{thm(1.12)} shows that Fuglede's conjecture holds for a special class of Moran measures $\mu_{\{b_k\},\{D_k\}}$.  In addition, this theorem can be used to prove that  $\mu_{\sigma}$ of Theorem \ref{thm(1.13)} is not a spectral measure if $\sigma\in\Pi_{l}$.
\begin{thm}\label{thm(1.12)}
Given  integers $p_{1}, p_{2},b_{1}\geq 2 $,  $t_{1}, t_{2}\geq 1$. Let  $b_{k}=p_{2} $ for all $k \geq 2$ and $\mu_{\{b_k\},\{D_k\}}$ be defined by \eqref{eq(2.8)}, where
 $$D_{k}= \begin{cases}\left\{0,1, \cdots, p_{1}-1\right\}t_{1}, & if~k=1 ;\\ \left\{0,1, \cdots, p_{2}-1\right\}t_{2} , & if~k\geq2.   \end{cases}$$
Then the following statements are equivalent.
\begin{enumerate}[(i)]
\item $t_{2}\mid t_{1}$.
 \item $\mu_{\{b_{k}\},\{D_k\}}$ is a spectral measure.
 \item $K(b_{k}, D_{k})$ is a translation tile, where $K(b_{k}, D_{k})$ is defined by \eqref{eq(0.0)}.
\end{enumerate}
\end{thm}
 \begin{rem}
By  the above Theorem \ref{thm(1.12)},  we can see  that the condition $p_{k+1}\mid(b_{k+1}t_{k })$ in Theorem \ref{thm(1.1)} is not sufficient  condition for $\mu_{\{b_k\},\{D_k\}}$ to be a spectral measure. But we also can  construct some examples such that the necessary condition in Theorem \ref{thm(1.1)} is also sufficient condition. For example, let $D_{2k-1}=\{0,1,2,3,4,5\}$, $D_{2k}=\{0,1\}3$ and $b_1=6$, $b_{2k+1}=2$, $b_{2k}=6$ for $k\geq 1$, then $p_{k+1}\mid(b_{k+1}t_{k})$ and $p_{k}\nmid t_{k+1}$ for $k\geq 1$, where $p_{2k-1}=6$, $p_{2k}=2$ and $t_{2k-1}=1$, $t_{2k}=3$. Moreover, we can rewrite $\mu_{\{b_k\},\{D_k\}}$ as follows form
$$
\begin{aligned}
\mu_{\{b_k\},\{D_k\}}&=\delta_{6^{-1}\{0,1,2,3,4,5\}}\ast\delta_{6^{-2}3\{0,1\}}
\ast\delta_{2^{-1}6^{-2}\{0,1,2,3,4,5\}}\ast\delta_{2^{-1}6^{-3}3\{0,1\}}\ast\cdots \\
&=\delta_{12^{-1}\{0,1,\cdots,11\}}\ast\delta_{12^{-2}\{0,1,\cdots,11\}}
\ast\delta_{12^{-3}\{0,1,\cdots,11\}}\ast\cdots\\
&:=\mu_{12,\{0,1,\cdots,11\}}.
\end{aligned}
$$
Then by using Theorem \ref{thm(1)}, one can conclude that $\mu_{\{b_k\},\{D_k\}}$ is a spectral measure.
\end{rem}

The paper is organized as follows. In Section 2,  we introduce some basic definitions and related  results that will be used in the proof of our main theorems. In Section 3, we study the structure of the spectrum of $\mu_{\{b_k\},\{D_k\}}$, and  then  prove Theorems \ref{thm(1.1)} and  \ref{thm(1.12)}. In Section 4,  
we focus on proving Theorem \ref{thm(1.13)}.  At the end of this paper, we  give some examples related to our main results.

\section{\bf Preliminaries \label{sect.2}}
In this section, we give some preliminary results and some  relevant concepts that will be used later.
Let the Moran measure  $\mu_{\{b_k\},\{D_k\}}$  be generated by an integer sequence $\{b_k\}_{k=1}^{\infty}$ with $|b_k|\geq 2$  and an integer sequence of digit sets $\{D_k\}_{k=1}^{\infty}$, where $D_{k}=\left\{0,1, \cdots, p_{k}-1\right\}t_{k}$ with $p_k\geq2$, $| t_{k}|\geq1$.
The Fourier transform $\widehat{\mu}_{\{b_{k}\},\{D_{k}\}}$ of $\mu_{\{b_k\},\{D_k\}}$ is given by
\begin{equation}\label{eq(2.1)}
\widehat{\mu}_{\{b_{k}\},\{D_{k}\}}(x)=\int e^{2\pi i x \xi }d\mu_{\{b_{k}\},\{D_{k}\}}(\xi)
=\prod\limits_{k=1}^{\infty}m_{D_{k}}\big(\frac{x}{b_1b_2\cdots b_k}\big) ,\quad x \in\mathbb{R},
\end{equation}
where
\begin{equation}\label{eq(2.2)}
  m_{D_{k}}(x):=\frac{1}{p_{k}}\sum_{j=0}^{p_{k}-1}e^{2\pi ijt_{k}x}, \quad x\in \mathbb{R}.
\end{equation}
It is easy to show that $m_{D_{k}}(x)$ is a $\frac{\mathbb{Z}}{t_{k}}$-periodic function.
 For a function $f(x) (x\in\mathbb{R})$, let $\mathcal{Z}(f)$ denote the zero sets of $f$, i.e., $\mathcal{Z}(f)=\big\{x\in\mathbb{R}:\, f(x)=0\big\}$. From \eqref{eq(2.1)}, we have
\begin{equation}\label{eq(2.3)}
\mathcal{Z}(\widehat{\mu}_{\{b_k\},\{D_k\}})=\left\{x \in \mathbb{R}:~ m_{D_{k}}\big(\frac{x}{b_1b_2\cdots b_k}\big)=0 \text{ for some}~k \in \mathbb{N}\right\}=\bigcup_{k=1}^{\infty}b_1b_2\cdots b_k\mathcal{Z}\big(m_{D_k}\big).
\end{equation}
By a direct calculation, it is easy to get that
\begin{equation}\label{eq(2.4)}
\mathcal{Z}(m_{D_k})=\frac{\mathbb{Z}\setminus p_{k}\mathbb{Z}}{p_{k}t_{k}}.
\end{equation}
Hence
\begin{equation}\label{eq(2.5)}
\mathcal{Z}(\widehat{\mu}_{\{b_k\},\{D_k\}})=\bigcup_{k=1}^{\infty}\left(b_1b_2 \cdots b_k \frac{\mathbb{Z}\setminus p_{k}\mathbb{Z}}{p_{k}t_{k}}\right) \subset \bigcup_{k=1}^{\infty}\frac{b_1 b_2\cdots b_k\mathbb{Z}}{p_{k}t_{k}}.
\end{equation}
For a countable discrete set $\Lambda\subset \mathbb{R}$, it is easy to see
that $E(\Lambda)=\{e^{2\pi i\lambda x}:\lambda\in \Lambda\}$  is an orthogonal family of $L^{2}(\mu_{\{b_k\},\{D_k\}})$ if and only if
$$
0=\langle e^{2\pi i \lambda_1 x},e^{2\pi i\lambda_2 x}\rangle_{L^2(\mu_{\{b_k\},\{D_k\}})}=\int e^{2\pi i (\lambda_1-\lambda_2) x}d\mu_{\{b_k\},\{D_k\}}=\widehat{\mu}_{\{b_k\},\{D_k\}}(\lambda_1-\lambda_2)
$$
for any $\lambda_1\neq\lambda_2\in\Lambda$. Therefore, the orthogonality of $E(\Lambda)$ is equivalent
to
\begin{equation}\label {eq(2.6)}
 (\Lambda-\Lambda)\setminus \{0\} \subset\mathcal{Z}(\widehat{\mu}_{\{b_k\},\{D_k\}}).
\end{equation}
\begin{defi}\label{prop(2.1)}
Let  $|b|\geq 2$  be an integer, and let $D, L\subset\mathbb{Z}$ be two finite digit
sets with the same cardinality. We say
that  $(b, D)$ forms an {\it  admissible pair}  (or $(b^{-1}D, L)$ forms a {\it compatible pair} or $(b, D, L)$ forms a {\it Hadamard triple}) if the matrix
$$
H=\frac{1}{\sqrt{\#D}}\left[  e^{2\pi i  dl/b }\right]_{d\in D, l\in L}
$$
is unitary, i.e., $H^*H=I$, where $H^*$ denotes the conjugate transposed matrix of $H$.
\end{defi}
The following properties of compatible pair can be found in \cite{DHL2019} or  be checked directly.
\begin{prop}\label{prop(2.2)}
Let  $|b|\geq 2$  be an integer, and let $D, L\subset\mathbb{Z}$ be two finite digit
sets with the same cardinality, then following statements are equivalent:
\begin{enumerate}[(i)]
 \item  $(b^{-1}D,  L)$ is a compatible pair.
 \item  $m_{D}(b^{-1}(l_{1}-l_{2}))=0$ for any $l_{1} \neq l_{2} \in L$.
 \item  $\sum_{l \in L} |m_{D } (\frac{l}{b}+x ) |^{2}=1$ for any $x\in\mathbb{R}$ .
\end{enumerate}
\end{prop}

For an  integer $|b| \geq2$ and an integer digit sets  $D=\{0,1, \ldots, p-1\}t$ with $p\geq2$, $|t|\geq1$, the following  proposition give an equivalent condition for $(b, D)$ to be an admissible pair.
\begin{prop}\label{prop(2.5)}Let $|b|\geq 2$  be an integer, and let integer digit sets $D=\{0,1, \ldots, p-1\}t$ with $p\geq2$, $|t|\geq1 $. Then $(b, D)$ is an admissible pair  if and only if
$ p\mid\frac{b }{\gcd(b , t )}$.
\end{prop}
\begin{proof}
Let $s=\gcd(b, t )$ and $b =b's $, $t =t's$.
We first prove the necessity.  Suppose that $(b , D )$ is an admissible pair. By Proposition \ref{prop(2.2)}, there exists $ L \subset\mathbb{Z}$ with $0\in L$ and $\#L=p$ such that
$$L \setminus\{0\}\subset (L -L )\setminus\{0\}\subset \frac{b (\mathbb{Z}\setminus p \mathbb{Z})}{tp}.$$
 This implies that there exist $k_{i}\in \mathbb{Z}$ for $i\in \{1,\cdots,p-1\}$ such that $$L=\left\{0,\frac{b(1+k_{1}p)}{tp},\cdots, \frac{b(p-1+k_{p-1}p)}{tp}\right\}\subset\mathbb{Z}.$$
Then $\frac{b(1+k_{1}p)}{tp}=\frac{b'(1+k_{1}p)}{t'p}\in\mathbb{Z}$.  According to $\gcd( 1+k_{1}p,p )=1$ and  $\gcd(b', t')=1$, we have $p\mid b'$, i.e., $p\mid\frac{b}{\gcd(b, t)}$.

Next, we prove the sufficiency. According to $p\mid\frac{b}{\gcd(b, t)}$, we have $\gcd(p, t')=1$. Hence, $t'\{0,1,\cdots,p-1\}\equiv \{0,1,\cdots,p-1\}({\rm mod}\ p)$. Let $L=\frac{bt'}{tp}\{0,1,\cdots,p-1\}$. It is easy to get $L\subset \mathbb{Z}$ and $(L-L)\setminus\{0\}\subset \frac{b(\mathbb{Z}\setminus p\mathbb{Z})}{tp}.$  It follows from Proposition \ref{prop(2.2)} that $(b^{-1}D,  L)$ is a compatible pair, and then  $(b , D )$ is admissible.
\end{proof}
In \cite[Lemma 4.2]{JP98}, Jorgensen and Pedersen  given a criterion that allows us to determine whether a countable set  $\Lambda$  is an orthonormal set or a spectrum of the measure $\mu$.
\begin{prop}\label{prop(2.3)}\cite{JP98}
Let $\mu$ be a Borel probability measure with compact support, and let $Q_{\mu,\Lambda}(x)=\sum_{\lambda \in \Lambda}|\widehat{\mu}(x+\lambda)|^{2}$ for  a countable set $\Lambda \subset \mathbb{R}.$ Then
\begin{enumerate}[(i)]
 \item $\Lambda$ is an orthonormal set of $\mu$ if and only if $Q_{\mu,\Lambda}(x) \leq 1$ for $x \in \mathbb{R}$.
 \item $\Lambda$ is a spectrum of $\mu$ if and only if $Q_{\mu,\Lambda}(x) \equiv 1$ for $x\in \mathbb{R}$.
\end{enumerate}
Moreover,  if $\Lambda$ is an orthonormal set, then $Q_{\mu,\Lambda}(x)$ is an entire function.
\end{prop}
The following lemma tells us that the spectrality of $\mu_{\{b_{k}\},\{D_{k}\}}$ is invariant under a linear transformation.
\begin{lem}\label{lem(2.3)}
Given a sequence   $\{b_k\}_{k=1}^{\infty}\subset \mathbb{R}$ with $|b_{k}|>1$ and a sequence of  digit sets $\{D_k\}_{k=1}^{\infty}$ with $D_k\subset\mathbb{R}$. Let $\tilde{D}_{k}=qD_{k}$ with $q\in\mathbb{R}\setminus\{0\}$. Then $\mu_{\{b_{k}\},\{D_{k}\}}$ is a spectral measure with spectrum $\Lambda$ if and only if  $\mu_{\{b_{k}\},\{\tilde{D}_{k}\}}$ is a spectral measure with spectrum $q^{-1}\Lambda$.
\end{lem}
\begin{proof}
For any $x\in\mathbb{R}$, we have
\begin{align*}
Q_{\mu_{\{b_{k}\},\{D_{k}\}}, \Lambda}(x)&=\sum_{\lambda \in \Lambda}|\widehat{\mu}_{\{b_{k}\},\{D_{k}\}}(x+\lambda)|^{2}=\sum_{\lambda \in \Lambda}|\widehat{\mu}_{\{b_{k}\},\{\tilde{D}_{k}\}}(q^{-1}x+q^{-1}\lambda)|^{2}
\\&=\sum_{\lambda' \in q^{-1}\Lambda}|\widehat{\mu}_{\{b_{k}\},\{\tilde{D}_{k}\}}(q^{-1}x+\lambda')|^{2}
=Q_{\mu_{\{b_{k}\},\{\tilde{D}_{k}\}}, q^{-1}\Lambda}(q^{-1}x).
\end{align*}
Hence, the conclusion follows from  Proposition \ref{prop(2.3)} $(ii)$.
\end{proof}
\begin{rem}\label{rem(2.3)}
For any $b\in\mathbb{R}\setminus\{0\}$. Let $\tilde{D}_{k}=\frac{ b_{1}}{b}D_{k}$, then
 \begin{align*}\mu_{\{b_{k}\},\{\tilde{D}_{k}\}}&=\delta_{ b_{1}^{-1}\frac{ b_{1}}{b}D_{1}}* \delta_{ b_{1}^{-1} b_{2}^{-1} \frac{ b_{1}}{b}D_{2}} * \cdots * \delta_{b_{1}^{-1} b_{2}^{-1} \cdots b_{k}^{-1}\frac{ b_{1}}{b}D_{k}} * \cdots
 \\&=\delta_{b^{-1}D_{1}}* \delta_{ b^{-1} b_{2}^{-1} D_{2}} * \cdots * \delta_{b^{-1}b_{2}^{-1} \cdots b_{k}^{-1} D_{k}} * \cdots.
 \end{align*}
According to Lemma \ref{lem(2.3)},  the spectrality of $\mu_{\{b_{k}\},\{D_{k}\}}$ is  equivalent to $\mu_{\{b_{k}\},\{\tilde{D}_{k}\}}$. This means that the spectrality of $\mu_{\{b_{k}\},\{D_{k}\}}$ is independent of $b_{1}$.
\end{rem}
The following proposition indicates that we can always assume that the  sequences $\{b_k\}_{k=1}^{\infty}$ and $\{t_k\}_{k=1}^{\infty}$ are positive in our main theorems.
\begin{prop}\label{prop(2.4)}
  Given integer sequences $\{b_k\}_{k=1}^{\infty}$, $\{p_k\}_{k=1}^{\infty}$ and $\{t_k\}_{k=1}^{\infty}$ with $p_{k} , |b_{k}| \geq 2$ and $|t_{k}| \geq 1$. Suppose that $\left\{p_{k}t_{k}\right\}_{k=1}^{\infty}$ is a bounded sequecne, then $\mu:=\delta_{ b_{1}^{-1} \left\{0,1, \cdots, p_{1}-1\right\} t_{1} }* \delta_{(b_{1} b_{2})^{-1} \left\{0,1, \cdots, p_{2}-1\right\}t_{2}}* \delta_{(b_{1}b_{2}b_{3})^{-1}\left\{0,1, \cdots, p_{3}-1\right\}t_{3}} * \cdots$ is a spectral measure if and
only if  $\nu:=\delta_{ |b_{1}|^{-1}\left\{0,1, \cdots, p_{1}-1\right\}\left|t_{1} \right| } * \delta_{ |b_{1}b_{2}|^{-1}\left\{0,1, \cdots, p_{2}-1\right\}\left|t_{2} \right| }  * \delta_{|b_{1}b_{2}b_{3}|^{-1}\left\{0,1, \cdots, p_{3}-1\right\}\left|t_{3} \right| } * \cdots$ is a spectral measure.

\end{prop}
\begin{proof}
For $k \geq 1$, let
$$
\gamma_{k}= \begin{cases}0, & if~~~b_{1} b_{2} \cdots b_{k}t_{k}>0; \\ -\left(b_{1} b_{2} \cdots b_{k}\right)^{-1}\left(p_{k}-1\right)t_{k}, & if~~~b_{1} b_{2} \cdots b_{k}t_{k}<0. \end{cases}
$$
 Since $\left|b_{k}\right| \geq 2$ and the sequence $\left\{p_{k}t_{k} \right\}_{k=1}^{\infty}$ is bounded, $\gamma:=\sum_{k=1}^{\infty}\gamma_{k}$ is a finite number.
For any $x\in\mathbb{R}$, it follows from \eqref{eq(2.1)} and  \eqref{eq(2.2)} that
\begin{equation}\label{eq(2.111)}
\begin{aligned}
\widehat{\nu}(x)&=\int e^{2\pi i x\xi }d\nu(\xi)
=\prod\limits_{k=1}^{\infty} \frac{1}{p_{k}}\sum_{j=0}^{p_{k}-1}e^{2\pi ix\mid\frac{ j t_{k}}{b_{1}b_{2}\cdots b_{k}}\mid}
\\&=e^{2\pi i \gamma x}\prod\limits_{k=1}^{\infty} \frac{1}{p_{k}}\sum_{j=0}^{p_{k}-1}e^{2\pi ix \frac{ j t_{k}}{b_{1}b_{2}\cdots b_{k}}}=e^{2\pi i\gamma x} \widehat{\mu}(x).
\end{aligned}
\end{equation}
Therefore, $\left|\widehat{\mu}(x)\right|=\left|\widehat{\nu}(x)\right|$ and the assertion follows easily from  Proposition \ref{prop(2.3)} $(ii)$.
\end{proof}

In \cite[Lemma 2.5]{DL21},  Deng et al. proved the following lemma which plays a key role in the proof of  Theorem \ref{thm(1.1)}.
\begin{lem}\label{lem(2.4)}\cite{DL21}
Let $p_{i, j}>0$ be positive numbers such that $\sum_{j=1}^{n} p_{i, j}=1(i=1,2, \cdots, m)$ and $\sum_{i=1}^{m} \max \left\{x_{i, j}: 1 \leq j \leq n\right\} \leq 1$ with $x_{i, j} \geq 0$. Then $\sum_{i=1}^{m} \sum_{j=1}^{n} p_{i, j} x_{i, j}=1$ if and only if  $\sum_{i=1}^{m} x_{i, 1}=1$  and $x_{i, 1}=x_{i, 2}=\cdots=x_{i,n}$ for $1 \leq i \leq m$.
\end{lem}

%

\section{\bf   Proof of  Theorems  \ref{thm(1.1)} and  \ref{thm(1.12)}} \label{sect.3}
In this section, we first investigate  the structure of the spectrum of $ \mu_{\{b_k\},\{D_k\}}$, and then prove Theorems  \ref{thm(1.1)} and  \ref{thm(1.12)}.

Given integer sequences $\{b_k\}_{k=1}^{\infty}$, $\{p_k\}_{k=1}^{\infty}$ and $\{t_k\}_{k=1}^{\infty}$ with $p_{k}$, $|b_{k}| \geq 2$  and $|t_{k}| \geq 1$. Let $\{D_k\}_{k=1}^{\infty}$ be an  integer sequence of digit sets, where $D_{k}=\left\{0,1, \cdots, p_{k}-1\right\}t_{k}$ and $\left\{p_{k}t_{k}\right\}_{k=1}^{\infty}$ is  a bounded sequence. It follows from \eqref{eq(2.82)} that
 $$\mu_{k}:=\delta_{b_{ 1}^{-1}D_{ 1}}\ast\delta_{b_{ 1}^{-1}b_{ 2}^{-1}D_{ 2}}
\ast\cdots\ast \delta_{b_{ 1}^{-1}b_{ 2}^{-1}\cdots b_{ k}^{-1} D_{k }}$$
converges weakly to Moran measure
\begin{equation}\label{eq(2.08)}\begin{aligned}
 \mu_{\{b_k\},\{D_k\}} =\delta_{b_{1}^{-1}D_1}\ast\delta_{b_{1}^{-1}b_{2}^{-1}D_2}
\ast\delta_{b_{1}^{-1}b_{2}^{-1} b_{3}^{-1}D_3}\ast\cdots=\mu_{k}\ast\nu_{>k}(b_{1}b_{2}\cdots b_{k}\cdot)\end{aligned}
\end{equation}
and $\mu_{\{b_k\},\{D_k\}}$  has a compact support, where
 \begin{equation}\label {eq(2.9)}
\nu_{>k}:=\delta_{b_{k+1}^{-1}D_{k+1}}\ast\delta_{b_{k+1}^{-1}b_{k+2}^{-1}D_{n+2}}
\ast\delta_{b_{k+1}^{-1}b_{k+2}^{-1}b_{k+3}^{-1}D_{k+3}}\ast\cdots.
\end{equation}
According to Proposition \ref{prop(2.4)}, we can always assume that $b_{k}\geq2$ and $t_{k}\geq1$ in the following study of the spectrality of $\mu_{\{b_k\},\{D_k\}}$.

Since the sequence $\left\{p_{k}t_{k}\right\}_{k=1}^{\infty}$ is bounded,  let $q \in \mathbb{N}$ be a common multiple of $\{p_{k}t_{k}\}_{k=1}^{\infty}$ and
$\tau_{k} =\frac{q}{p_{k}t_{k}}$ for $k\geq1$. Assume that $\Lambda$ is a spectrum of $\mu_{\{b_k\},\{D_k\}} $ with $0 \in \Lambda$.
According to \eqref{eq(2.5)}, we have $\frac{q}{b_{1}}\Lambda\subset \mathbb{Z}$. This means that $\frac{1}{b_{1}}\Lambda\subset\frac{\{0,1,\cdots,q-1\}}{q}+  \mathbb{Z}$. Hence
$$
\frac{1}{b_{1}} \Lambda=\bigcup_{n=0}^{q-1}\left(\frac{n}{q}+ \Lambda_{n}\right),
$$
where $ \Lambda_{n} =\mathbb{Z} \cap(\frac{\Lambda}{b_{1}} - \frac{n}{q}) $   and  $\frac{n}{q}+\Lambda_{n}$ is empty when $\Lambda_{n}=\emptyset$.
For any $n\in\{0,1,\cdots,q-1\}$, there exist  unique $i\in\{0,1,\cdots,\tau_{1}-1\}, j\in\{0,1,\cdots,p_{1}-1\}$ and $l\in\{0,1,\cdots,t_{1}-1\}$ such that $n=i+\tau_{1} j+\tau_{1} p_{1}l$.
Thus we have the following decomposition
\begin{equation}\label{eq(3.1)}
\frac{1}{b_{1}} \Lambda=\bigcup_{n=0}^{q-1}\left(\frac{n}{q}+ \Lambda_{n}\right)=\bigcup_{i=0}^{\tau_{1}-1} \bigcup_{j=0}^{p_{1}-1} \bigcup_{l=0}^{t_{1}-1}\left(\frac{i+\tau_{1}j+\tau_{1} p_{1}l}{q}+\Lambda_{i+\tau_{1} j+\tau_{1} p_{1}l}\right).
\end{equation}

Under the assumption that $\mu_{\{b_k\},\{D_k\}} $  is a spectral measure, the following proposition characterizes  the structure of the spectra of $\nu_{>1}$.
\begin{prop}\label{prop(3.3)}
Let $\mu_{\{b_k\},\{D_k\}}$ and $\nu_{>1}$  be defined by \eqref{eq(2.08)} and \eqref{eq(2.9)} respectively. Suppose $\Lambda$ is a spectrum of $\mu_{\{b_k\},\{D_k\}}$ with $0 \in \Lambda$.  For any group $\{j_i:0\leq i\leq\tau_{1}-1\}\subset \{0,1,\cdots,p_{1}-1\}$,  write $$\Gamma=\bigcup_{i=0}^{\tau_{1}-1}\bigcup_{l=0}^{t_{1}-1}\left(\frac{i+\tau_{1}j_{i}+\tau_{1} p_{1}l}{q}+\Lambda_{i+\tau_{1} j_{i}+\tau_{1} p_{1}l}\right).$$
 Then $\Gamma$  is a spectrum of $\nu_{>1}$ if $\Gamma\neq\emptyset$.
\end{prop}
\begin{proof}We will divide the  proof into the following two steps.

\textbf{Step 1.} $\Gamma $ is an orthogonal set of $\nu_{>1}$.

The conclusion obviously holds if $\Gamma $ contains only  one element.
Next, we consider the case that $\Gamma $ contains  at least two elements. For any distinct $\lambda_1, \lambda_2\in \Gamma $, we can write $$\lambda_1= \frac{i_{1}+\tau_{1}j_{i_{1}}+\tau_{1} p_{1}l_{1}}{q}+z_{1} \quad  \text{and } \quad \lambda_2=\frac{i_{2}+\tau_{1} j_{i_{2}}+\tau_{1} p_{1}l_{2}}{q}+z_{2}$$
for some $i_{1},  i_{2} \in\{0,1,\cdots,\tau_{1}-1\}$, $l_{1}, l_{2}\in\{0,1,\cdots,t_{1}-1\}$, $z_1\in\Lambda_{i_{1}+\tau_{1} j_{i_{1}}+\tau_{1} p_{1}l_{1}} $ and  $z_2\in\Lambda_{i_{2}+\tau_{1}j_{i_{2}}+\tau_{1} p_{1}l_{2}}$.
 Since $\Lambda$ is a spectrum of $\mu_{\{b_k\},\{D_k\}}$ with $0 \in \Lambda$, it follows from $\eqref{eq(3.1)}$ that $$\Lambda= \bigcup_{i=0}^{\tau_{1}-1} \bigcup_{j=0}^{p_{1}-1} \bigcup_{l=0}^{t_{1}-1}\left(\frac{b_{1}(i+\tau_{1}j+\tau_{1}p_{1}l)}{q}+b_{1}\Lambda_{i+\tau_{1}j+\tau_{1}p_{1}l}\right)$$  and $b_1\Gamma \subset \Lambda$. Hence the orthogonality  of $\Lambda$ implies that
 $$
 0=\widehat{\mu}_{\{b_k\},\{D_k\}} (b_{1}(\lambda_1-\lambda_2) )=m_{D_1} (\lambda_1-\lambda_2)  \widehat{\nu}_{>1} (\lambda_1-\lambda_2).
 $$
According to the fact that $m_{D_1}$ is a $\frac{\mathbb{Z}}{t_{1}}$-periodic function, we have
$$
\begin{aligned}
 m_{D_1} (\lambda_1-\lambda_2 )& = m_{D_1}\left( \frac{i_{1}-i_{2}+\tau_{1} (j_{i_{1}}-j_{i_{2}})+\tau_{1} p_{1}(l_{1}-l_{2})}{q}   \right)
 \\&= m_{D_1}\left( \frac{i_{1}-i_{2}+\tau_{1} (j_{i_{1}}-j_{i_{2}}) }{q}  \right).
\end{aligned}$$

If $i_{1}= i_{2}$, then $j_{i_{1}}=j_{i_{2}}$ and $ m_{D_1} (\lambda_1-\lambda_2)=m_{D_1}(0)=1$. This implies that $\widehat{\nu}_{>1} ( \lambda_1-\lambda_2 )=0.$

If $i_{1}\neq i_{2}$,  then  $i_{1}-i_{2}\in\{1-\tau_{1},2-\tau_{1},\cdots,-1,1,\cdots,\tau_{1}-2, \tau_{1}-1\}$. This yields that $i_{1}-i_{2}\notin\tau_{1} \mathbb{Z}$. We claim that
$\frac{i_{1}-i_{2}+\tau_{1}(j_{i_{1}}-j_{i_{2}}) }{q} \notin   \frac{ \mathbb{Z}\setminus p_{1} \mathbb{Z}}{p_{1} t_{1}}=\mathcal{Z}(m_{D_1}).$ If not , there exists $k\in\mathbb{Z}\setminus p_{1} \mathbb{Z}$ such that  $\frac{i_{1}-i_{2}+\tau_{1} (j_{i_{1}}-j_{i_{2}}) }{q}=\frac{k}{p_{1}t_{1}}.$ Note that $q=p_{1}t_{1}\tau_{1}$, we derive that $i_{1}-i_{2}=\tau_{1}(k-j_{i_{1}}+j_{i_{2}})$, which contradicts with $i_{1}-i_{2}\notin\tau_{1} \mathbb{Z}$.  So the claim follows.
Hence $ m_{D_1} ( \lambda_1-\lambda_2 ) = m_{D_1}\left( \frac{i_{1}-i_{2}+\tau_{1} (j_{i_{1}}-j_{i_{2}}) }{q}  \right) \neq 0. $
This also implies that $\widehat{\nu}_{>1}\left( \lambda_1-\lambda_2 \right)=0.$

Therefore, $\Gamma $ is an orthogonal set of $\nu_{>1}$.

 \textbf{Step 2.  } $Q_{  \nu_{>1},\Gamma}\left(  x\right)=\sum_{\gamma\in\Gamma }\left|\widehat{\nu }_{>1}\left(x+\gamma\right)\right|^{2}\equiv1$.

Write $\tilde{\Lambda}_{ijl}:=\Lambda_{i+\tau_{1}j+\tau_{1}p_{1}l}$.
By Proposition \ref{prop(2.3)} $(ii)$ and \eqref{eq(3.1)}, it follows that
$$
1=Q_{\mu_{\{b_k\},\{D_k\}},\Lambda}\left(b_{1} x\right)=\sum_{i=0}^{\tau_{1}-1} \sum_{j=0}^{p_{1}-1} \sum_{l=0}^{t_{1}-1}\sum_{\lambda \in(i+\tau_{1}j+\tau_{1}p_{1}l)+q\tilde{\Lambda}_{ijl}}\left|\widehat{\mu}_{\{b_k\},\{D_k\}}\left(\frac{b_{1}}{q} \lambda+b_{1} x\right)\right|^{2}
$$
for any $x \in \mathbb{R}$, where $\sum_{\lambda \in i+\tau_{1}j+\tau_{1} p_{1}l +q \tilde{\Lambda}_{ijl}}\left|\widehat{\mu}_{\{b_k\},\{D_k\}}\left(\frac{b_{1}}{q} \lambda+b_{1} x\right)\right|^{2}=0$ if $\tilde{\Lambda}_{ijl}=\emptyset$.

Using the fact that $\widehat{\mu}_{\{b_k\},\{D_k\}}(x)=m_{D_1}(b_{1}^{-1} x) \widehat{\nu}_{>1}(b_{1}^{-1} x)$ and  $m_{D_{1}}(x)$ is a $\frac{\mathbb{Z}}{t_{1}}$-periodic function, we have
\begin{footnotesize}
\begin{equation}\label{eq(3.222)}
\begin{aligned}
1&=\sum_{i=0}^{\tau_{1}-1} \sum_{j=0}^{p_{1}-1}\sum_{l=0}^{t_{1}-1} \sum_{\tilde{\lambda} \in \tilde{\Lambda}_{ijl}}\left| m_{D_1}\left(\frac{1}{q}(i+\tau_{1}j+\tau_{1} p_{1}l)+\tilde{\lambda}+x\right)\right|^{2}\left|\widehat{\nu}_{>1}\left(\frac{1}{q}(i+\tau_{1}j+\tau_{1} p_{1}l)+\tilde{\lambda}+x\right)\right|^{2}\\
&=\sum_{i=0}^{\tau_{1}-1}\sum_{j=0}^{p_{1}-1}\left| m_{D_1}\left(\frac{1}{q}(i+\tau_{1}j)+x\right)\right|^{2} \sum_{l=0}^{t_{1}-1}\sum_{\tilde{\lambda} \in\tilde{ \Lambda}_{ijl}}\left|\widehat{\nu}_{>1}\left(\frac{1}{q}(i+\tau_{1}j+\tau_{1} p_{1}l)+\tilde{\lambda}+x\right)\right|^{2}
\end{aligned}
\end{equation}
\end{footnotesize}
for any $x \in \mathbb{R}$.

Let $x\in\mathbb{R}\setminus\mathbb{Q}$.
Write $x_{i, j}=\sum_{l=0}^{t_{1}-1}\sum_{\tilde{\lambda}\in \tilde{\Lambda}_{ijl }}\left|\widehat{\nu}_{>1}\left(\frac{1}{q}(i+\tau_{1}j +\tau_{1} p_{1}l)+\tilde{\lambda}+x\right)\right|^{2}$ and $p_{i, j}=\left|m_{D_1}\left(\frac{1}{q}(i+ \tau_{1}j)+x\right)\right|^{2}.$
Then one may derive from \eqref{eq(2.4)} that $p_{i, j}>0$, and \eqref{eq(3.222)}  can be rewrite as
\begin{equation}\label{eq(3.223)}
\sum_{i=0}^{\tau_{1}-1}\sum_{j=0}^{p_{1}-1}p_{i, j}x_{i, j}=1.
\end{equation}
Define $L_{i}=\left\{i+\tau_{1}j: 0 \leq j \leq p_{1}-1\right\}$. Since $(q^{-1} D_{1}, L_{i})$ is a compatible pair for any $i\in\{0,1, \cdots, \tau_{1}-1\}$,  it follows from Proposition \ref{prop(2.2)} that
\begin{equation}\label{eq(3.224)}
\sum_{j=0}^{p_{1}-1} p_{i, j}=\sum_{j=0}^{p_{1}-1}  \left|m_{D_1}\left(\frac{1}{q}(i+ \tau_{1}j)+x\right)\right|^{2}=1
\end{equation}
for any  $i\in\{0,1, \cdots, \tau_{1}-1\}$. Since $\Gamma=\bigcup_{i=0}^{\tau_{1}-1}\bigcup_{l=0}^{t_{1}-1}\left(\frac{i+\tau_{1}j_{i}+\tau_{1} p_{1}l}{q}+\Lambda_{i+\tau_{1} j_{i}+\tau_{1} p_{1}l}\right)$ is an orthogonal set of $\nu_{>1}$ for any group $\{j_i:0\leq i\leq\tau_{1}-1\}\subset \{0,1,\cdots,p_{1}-1\}$, we conclude from Proposition \ref{prop(2.3)} $(i)$  that
\begin{equation}\label{eq(3.225)}\sum_{i=0}^{\tau_{1}-1} \max \left\{x_{i, 0}, x_{i, 1}, \cdots, x_{i, p_{1}-1}\right\} \leq 1.
\end{equation}
According to \eqref{eq(3.223)},  \eqref{eq(3.224)},  \eqref{eq(3.225)} and Lemma \ref{lem(2.4)}, it can be deduced that
\begin{equation}\label{eq(3.21)}
 \sum_{i=0}^{\tau_{1}-1}\sum_{l=0}^{t_{1}-1}\sum_{\tilde{\lambda}\in \tilde{\Lambda}_{ijl}}\left|\widehat{\nu}_{>1}\left(\frac{1}{q}(i+\tau_{1} j+\tau_{1} p_{1}l)+\tilde{\lambda}+x\right)\right|^{2}=1, \quad (j=0,1, \cdots, p_{1}-1)
 \end{equation}
and
 \begin{footnotesize}
 \begin{equation}\label{eq(3.22)}
 \begin{aligned}
 \begin{aligned}\sum_{l=0}^{t_{1}-1}\sum_{\tilde{\lambda}\in \tilde{\Lambda}_{i0l }}\left|\widehat{\nu}_{>1}\left(\frac{1}{q}(i+ \tau_{1} p_{1}l) +\tilde{\lambda}+x\right)\right|^{2}&=  \sum_{l=0}^{t_{1}-1}\sum_{ \tilde{\lambda} \in\tilde{ \Lambda}_{i1l }}\left|\widehat{\nu}_{>1}\left(\frac{1}{q}(i+\tau_{1}  +\tau_{1}p_{1}l) +\tilde{\lambda}+x\right)\right|^{2}\\&=\cdots\\&= \sum_{l=0}^{t_{1}-1}\sum_{\tilde{\lambda}\in \tilde{\Lambda}_{i (p_{1}-1)l }}\left|\widehat{\nu}_{>1}\left(\frac{1}{q}\left(i+\tau_{1}( p_{1}-1)+\tau_{1} p_{1}l\right) +\tilde{\lambda}+x\right)\right|^{2}\end{aligned}
 \end{aligned}
 \end{equation}
 \end{footnotesize}
 for any $i\in\{0,1, \cdots, \tau_{1}-1\}$.
From the above equations \eqref{eq(3.21)} and \eqref{eq(3.22)}, we can know that
 \begin{equation}\label{eq(3.221)}
Q_{\nu_{>1}, \Gamma}\left(  x\right)=\sum_{\gamma\in\Gamma }\left|\widehat{\nu }_{>1}\left(x+\gamma\right)\right|^{2} =1
 \end{equation}
 for any $x\in\mathbb{R}\setminus\mathbb{Q}$.
Since $Q_{\nu_{>1},\Gamma}(x)$ is an entire function, the above equation \eqref{eq(3.221)}  holds  for any $x \in \mathbb{R}$. It  follows from Proposition \ref{prop(2.3)} $(ii)$  that $\Gamma$ is a spectrum of $ \nu_{>1}.$
\end{proof}
\begin{rem}\label{rem(3.4)}
 Suppose $\Lambda= \bigcup_{i=0}^{\tau_{1}-1} \bigcup_{j=0}^{p_{1}-1} \bigcup_{l=0}^{t_{1}-1}\left(\frac{b_{1}(i+\tau_{1}j+\tau_{1} p_{1}l)}{q}+b_{1}\Lambda_{i+\tau_{1}j+\tau_{1} p_{1}l}\right)$ is a spectrum of $\mu_{\{b_k\},\{D_k\}}$ with $0 \in \Lambda$. Then we can conclude from \eqref{eq(3.22)}  that for any $i\in \{0,1,\cdots,\tau_{1}-1\}$, one of the following two statements holds:
 \begin{enumerate}[(i)]
 \item    $\bigcup_{l=0}^{t_{1}-1}\Lambda_{i+\tau_{1} j+\tau_{1}p_{1}l}\neq\emptyset$ for all $j\in \{0,1,2,\cdots,p_{1}-1\}$.
 \item    $\bigcup_{l=0}^{t_{1}-1}\Lambda_{i+\tau_{1}j+\tau_{1} p_{1}l}=\emptyset$ for all $j\in \{0,1,2,\cdots,p_{1}-1\}$.
 \end{enumerate}
\end{rem}
\begin{prop}\label{prop(3.4)}
Let $\mu_{\{b_k\},\{D_k\}}$ be defined by \eqref{eq(2.08)}. Suppose that $0\in \Lambda$ is a spectrum of $\mu_{\{b_k\},\{D_k\}}$, then for any $j\in \{1,2,\cdots,p_{1}-1\}$, there exist
$l_{j}\in \{0,1,\cdots,t_{1}-1\}$ and $z_{j}\in \mathbb{Z}$ such that $\alpha_{j}:=\frac{b_{1}(j+ p_{1}l_{j})}{p_{1}t_{1}}+b_{1}z_{j} \in\Lambda$.
\end{prop}
\begin{proof}
It follows from \eqref{eq(3.1)} that $
\Lambda=\bigcup_{i=0}^{\tau_{1}-1} \bigcup_{j=0}^{p_{1}-1} \bigcup_{l=0}^{t_{1}-1}\left(\frac{b_{1}(i+\tau_{1}j+\tau_{1} p_{1}l)}{q}+b_{1}\Lambda_{i+\tau_{1}j+\tau_{1} p_{1}l}\right)$. Since $0\in\Lambda$,
then $0\in\bigcup_{l=0}^{t_{1}-1}\Lambda_{\tau_{1}p_{1}l}\neq\emptyset$. By Remark \ref{rem(3.4)}, we have $\bigcup_{l=0}^{t_{1}-1}\Lambda_{\tau_{1}j+\tau_{1}p_{1} l}\neq\emptyset$ for any $j\in \{1,2,\cdots,p_{1}-1\}$. Therefore,  for any $j\in \{1,2,\cdots,p_{1}-1\}$, there exist
$l_{j}\in \{0,1,\cdots,t_{1}-1\}$ and $z_{j}\in \Lambda_{\tau_{1}j +\tau_{1} p_{1}l_{j}} \subset\mathbb{Z}$ such that $\alpha_{j}:=\frac{b_{1}( j+ p_{1}l_{j})}{p_{1}t_{1}}+b_{1}z_{j} \in\Lambda$.
\end{proof}
\begin{prop}\label{prop(3.3.2)}
Let $\mu_{\{b_k\},\{D_k\}}$ and $\nu_{>k}$ be defined by \eqref{eq(2.08)} and \eqref{eq(2.9)}, respectively. Suppose $\mu_{\{b_k\},\{D_k\}}$  is a spectral measure, then $\nu_{>k}$ is a spectral measure for all $k\geq1$.
\end{prop}
\begin{proof}Since $\mu_{\{b_k\},\{D_k\}}$  is a spectral measure,  it follows from Proposition \ref{prop(3.3)} that  $\nu_{>1}$ is a spectral measure. Applying this conclusion to $\nu_{>1}$, then $\nu_{>2}$  is also a spectral measure. Repeat this operations, we can get that $\nu_{>k }$ is a spectral measure for all $k \geq 1 $.
\end{proof}

\begin{prop}\label{prop(3.4.2)}

Let $\mu_{\{b_k\},\{D_k\}}$  be defined by \eqref{eq(2.08)}.  Suppose that $\mu_{\{b_k\},\{D_k\}}$ is a spectral measure and $p_{k_{0} }\nmid t_{k_{0}+1}$ for some $k_{0}\geq1$, then $p_{k_{0}+1}\mid(b_{k_{0}+1}t_{k_{0}})$.
\end{prop}
\begin{proof}
Without loss of generality, we can assume that $k_{0}=1$ by Proposition \ref{prop(3.3.2)}.
Let $0\in\Lambda$ be a spectrum of $\mu_{\{b_k\},\{D_k\}}$,
then $$
\Lambda=\bigcup_{i=0}^{\tau_{1}-1} \bigcup_{j=0}^{p_{1}-1} \bigcup_{l=0}^{t_{1}-1}\left(\frac{b_{1}(i+\tau_{1}j+\tau_{1} p_{1}l)}{q}+b_{1}\Lambda_{i+\tau_{1}j+\tau_{1} p_{1}l}\right).
$$
By Proposition \ref{prop(3.3)}, for any group $\{j_i:0\leq i\leq\tau_{1}-1\}\subset \{0,1,\cdots,p_{1}-1\}$, the set
\begin{equation}\label{eq(3.110)}
\Gamma_{j_{0}, j_{1}, \cdots, j_{\tau_{1}-1}}:=\bigcup_{i=0}^{\tau_{1}-1}\bigcup_{l=0}^{t_{1}-1}\left(\frac{1}{q}\left(i+\tau_{1}j_{i}+\tau_{1} p_{1}l \right)+\Lambda_{i+\tau_{1}j_{i}+\tau_{1}p_{1}l}\right)
\end{equation}
is a spectrum of $\nu_{>1}$ if $\Gamma_{j_{0}, j_{1}, \cdots, j_{\tau_{1}-1}}\neq\emptyset$. Since $0\in\Lambda$, it is easy to see that $0\in\Gamma_{0, j_{1}, \cdots, j_{\tau_{1}-1}}$ (i.e, $j_0=0$) and it is a spectrum of $\nu_{>1}$.

In order to get the conclusion $p_{2}\mid(b_{2}t_{1})$, we first prove the following claim.

\noindent\emph{{\bf Claim 1.} If $ p_{1} \nmid t_{2}$, then for any $i\in\{1,2,\cdots,\tau_{1}-1\}$, there exists $j_{i}\in\{0,1,\cdots,p_{1}-1\}$ such that  $\frac{t_{1} b_{2}}{p_{2}}\neq t_{2}\frac{ i+\tau_{1} j_i}{ p_{1}\tau_{1}}({\rm mod} 1)$.}
\begin{proof}[Proof of Claim 1]
Suppose, on the contrary, that there exists  $i_{0}\in\{1,2,\cdots,\tau_{1}-1\}$ such that  $\frac{t_{1} b_{2}}{p_{2}}=t_{2}\frac{ i_{0}+\tau_{1}j}{ p_{1}\tau_{1}}({\rm mod} 1)$ for any  $j\in\{0,1,\cdots,p_{1}-1\}$. Since $p_{1}\geq2$, then $j$ can be taken $0$ and $1$ at least. This means that
$$\frac{t_{1} b_{2}}{p_{2}}= t_{2}\frac{ i_{0}  }{ p_{1}\tau_{1}}+n_{0}
\quad\quad   \text{and}\quad\quad
 \frac{t_{1} b_{2}}{p_{2}}= t_{2}\frac{ i_{0}+\tau_{1} }{ p_{1}\tau_{1}}+n_{1} $$
for some $n_{0}, n_{1}\in \mathbb{Z}$. The above two equations imply that $\frac{t_{2}  }{ p_{1}}=n_{0}-n_{1}\in  \mathbb{Z}$, which contradicts with $ p_{1} \nmid t_{2}$. Hence the claim follows.\end{proof}
According to Claim 1, for any  $i\in\{1,2,\cdots,\tau_{1}-1\}$, we choose a $\tilde{j}_{i}\in\{0,1,\cdots,p_{1}-1\}$ such that $\frac{t_{1} b_{2}}{p_{2}}\neq t_{2}\frac{ i+\tau_{1}\tilde{j_i}}{ p_{1}\tau_{1}}({\rm mod} 1)$. Let $\Gamma_{0,\tilde{j}_{1}, \cdots, \tilde{j}_{\tau_{1}-1}}$ be defined as \eqref{eq(3.110)}, then $\Gamma_{0,\tilde{j}_{1}, \cdots, \tilde{j}_{\tau_{1}-1}}$ is a spectrum of $\nu_{>1}$.
Applying Proposition \ref{prop(3.4)} to $\nu_{>1}$ and the spectrum $\Gamma_{0,\tilde{j}_{1}, \cdots, \tilde{j}_{\tau_{1}-1}}$, we conclude that there exist
$l_{j}\in \{0,1,\cdots,t_{2}-1\}$ and $z_{j}\in\mathbb{Z}$ such that
 $\alpha_{j}=\frac{b_{2}(j + p_{2}l_{j})}{p_{2}t_{2}}+b_{2}z_{j} \in\Gamma_{0, \tilde{j}_{1}, \cdots, \tilde{j}_{\tau_{1}-1}}$ for any $j\in \{1,2,\cdots,p_{2}-1\}$.
 Hence for $\alpha_{1}$, there exist  $i\in\left\{0,1, \cdots, \tau_{1}-1\right\}$,
$\tilde{l}\in\left\{0, 1, \cdots, t_{1}-1\right\}$  and $\tilde{z}_1\in \mathbb{Z}$ such that
$$
\alpha_{1}=\frac{b_{2}\left(1+p_{2}l_{1}\right)}{p_{2} t_{2}}+b_{2} z_{1}=\frac{i+\tau_{1}\tilde{j}_i+\tau_{1} p_{1}\tilde{l}}{q}+\tilde{z}_1\in \Gamma_{0, \tilde{j}_{1}, \cdots, \tilde{j}_{\tau_{1}-1}},
$$ where $\tilde{j}_0=0$.
Multiplying both sides of the above equation by $t_1t_2$, a simple calculation shows that
\begin{equation}\label{eq(3.120)}
\frac{t_{1} b_{2}}{p_{2}}= t_{2}\frac{ i+\tau_{1}\tilde{j}_i}{ p_{1}\tau_{1}}+ t_{2}\tilde{l} - t_{1}b_{2}l_{1}+ t_{2} t_{1}(\tilde{z}_1-z_1b_2).
\end{equation}
Since $\frac{t_{1} b_{2}}{p_{2}}\not\equiv t_{2}\frac{ i+\tau_{1}\tilde{j_i}}{ p_{1}\tau_{1}}({\rm mod} 1)$ for all $i\in\{1,2,\cdots,\tau_{1}-1\}$,  \eqref{eq(3.120)} holds only for $i=0$. Note that $\tilde{j}_0=0$, \eqref{eq(3.120)} implies that
$\frac{t_{1} b_{2}}{p_{2}}=t_{2}\tilde{l} - t_{1}b_{2}l_{1}+ t_{2} t_{1}(\tilde{z}_1-z_1b_2)\in \mathbb{Z}$. Hence $p_{2}\mid(b_{2}t_{1})$ and the proof is completed.
\end{proof}
Now, we can use the above results to prove Theorems  \ref{thm(1.1)} and  \ref{thm(1.12)}.
\begin{proof}[\bf Proof of Theorem \ref{thm(1.1)}]The proof follows directly from Propositions \ref{prop(3.3.2)} and \ref{prop(3.4.2)}.
\end{proof}

\begin{proof}[\bf Proof of Theorem \ref{thm(1.12)}]We first prove  $``(i)\Longleftrightarrow(ii) "$,  and then $``(i)\Longleftrightarrow(iii) "$.

 $``(i)\Longrightarrow(ii)"$. Let $t=\frac{t_{1}}{t_{2}}$ and $b = p_{1}t$. According to Lemma  \ref{lem(2.3)}, we can obtain that  $\mu_{\{b_{k}\},\{D_{k}\}}$ is a spectral measure  if and only if
 $$\begin{aligned}
\mu': =\delta_{b^{-1}\left\{0,1, \cdots,p_{1}-1\right\}t }\ast\delta_{b ^{-1}p_{2}^{-1}\left\{0,1, \cdots,p_{2}-1\right\} }
\ast\delta_{b^{-1}p_{2}^{-2}\left\{0,1, \cdots,p_{2}-1\right\} }\ast\delta_{b^{-1}p_{2}^{-3}\left\{0,1, \cdots,p_{2}-1\right\} }\ast \cdots
 \end{aligned}$$
is a spectral measure.
 Let $\Lambda= \{0,1,\cdots,p_{1}-1\}+b \mathbb{Z}$ and $D=\{0,1, \cdots,p_{1}-1 \}t$.
 Since $m_{D }(x)$ is a $\mathbb{Z}$-periodic function, we have
$$
Q_{\mu' ,\Lambda} (x) =\sum_{i=0}^{p_{1}-1} \sum_{\lambda'\in \mathbb{Z}}|\widehat{\mu'}  ( i +\lambda' +x ) |^{2}=\sum_{i=0}^{p_{1}-1}| m_{D} (\frac{i+x}{b}) |^{2}  \sum_{\lambda'\in \mathbb{Z}}  |\widehat{\nu'}_{>1} ( \frac{i +x}{b}+\lambda' ) |^{2},
$$
where $\nu'_{>1}=\delta_{p_{2}^{-1}\left\{0,1, \cdots,p_{2}-1\right\} }
\ast\delta_{p_{2}^{-2}\left\{0,1, \cdots,p_{2}-1\right\} }\ast\delta_{ p_{2}^{-3}\left\{0,1, \cdots,p_{2}-1\right\} }\ast \cdots$.
 By Proposition \ref{prop(2.2)}, it is  easy to show that $\sum_{i=0}^{p_{1}-1}| m_{D} (\frac{i+x}{b} ) |^{2}=1$.  It is well known that $\mathbb{Z}$  is a spectrum of  $\nu'_{>1}$. Then Proposition \ref{prop(2.3)} $(ii)$ implies that  $$Q_{\mu' ,\Lambda} (x)=\sum_{i=0}^{p_{1}-1}| m_{D} (\frac{i+x}{b } ) |^{2}  \sum_{\lambda'\in \mathbb{Z}} |\widehat{\nu}_{>1} (\frac{i +x}{b}+ \lambda') |^{2}=\sum_{i=0}^{p_{1}-1}| m_{D} (\frac{i+x}{b } ) |^{2}=1.$$ Therefore $ \mu '$ is  a spectral measure, so is
$\mu_{\{b_{k}\},\{D_{k}\}}$.

 $``(ii)\Longrightarrow(i)"$. Let  $q \in \mathbb{N}$ be a common multiple of $p_{1}t_{1}$ and $p_{2}t_{2}$, and let
$\tau_{1} =\frac{q}{p_{1}t_{1}}$. Suppose $\Lambda$ be a spectrum of  $\mu_{\{b_{k}\},\{D_k\}}$  with $0 \in \Lambda$. It follows from  Proposition \ref{prop(3.3)}
that for any group $\{j_i: 1\leq i\leq\tau_{1}-1\}\subset \{0,1,\cdots,p_{1}-1\}$, the nonempty set
$$
\begin{aligned}
 0\in \Gamma_{0, j_{1}, \cdots, j_{\tau_{ 1}-1}}:&=\bigcup_{i=0}^{\tau_{ 1}-1}\bigcup_{l=0}^{t_{ 1}-1}\left(\frac{1}{\tau_{ 1}p_{ 1}t_{ 1}}\left(i+\tau_{ 1}j_{i}+\tau_{ 1}p_{ 1}l  \right)+\Lambda_{i+\tau_{ 1}j_{i}+\tau_{ 1}p_{ 1}l}\right)
 \\&= \bigcup_{l=0}^{t_{ 1}-1} (\frac{l}{ t_{ 1}} +\Lambda_{ \tau_{ 1}p_{ 1}l} )\bigcup\left(\bigcup_{i=1}^{\tau_{ 1}-1}\bigcup_{l=0}^{t_{ 1}-1}\left (\frac{1}{\tau_{ 1}p_{ 1}t_{ 1}} (i+\tau_{ 1}j_{i}+\tau_{ 1}p_{ 1}l )+\Lambda_{i+\tau_{ 1}j_{i}+\tau_{ 1}p_{ 1}l}\right)\right)
 \end{aligned}
 $$
 is a spectrum of
 $\nu_{>1}=\delta_{p_{2}^{-1} D_{2}}  \ast \delta_{p_{2}^{-2} D_{2}} \ast \delta_{p_{2}^{-3}D_{2}} \ast \cdots$.

It is well known that $ \frac{\mathbb{Z}}{t_{2}}$  is the unique spectrum of  $\nu_{>1}$ containing zero.
Hence, we have $\bigcup_{l=0}^{t_{ 1}-1}\Lambda_{i+\tau_{ 1}j_{i}+\tau_{ 1}p_{ 1}l }=\emptyset$ for $i\in\{1,2,\cdots,\tau_{1}-1\}$. Otherwise, $j_{i}$ can be arbitrarily chosen from $\left\{0, 1,\cdots, p_{ 1}-1\right\}$ for $i\in\{1,2,\cdots,\tau_{1}-1\}$, which contradicts the fact that  $\nu_{>1}$  has the unique spectrum containing zero.  This implies that
$$\bigcup_{l=0}^{t_{ 1}-1}(\frac{l}{t_{ 1}}+\Lambda_{ \tau_{ 1}p_{ 1}l})= \frac{\mathbb{Z}}{t_{2}}.$$
Then there exist $\tilde{l}\in\{0,1,\cdots,t_{1}-1\} $ and $z_{0}\in\Lambda_{\tau_{ 1}p_{ 1}j}\subset\mathbb{Z}$ such that $\frac{\tilde{l}}{t_{ 1}}+z_{0}= \frac{1}{t_{2}}$. Hence  $\tilde{l}+z_{0}t_{ 1}= \frac{t_{ 1}}{t_{2}} $. One may infer that $t_{2}\mid t_{ 1}$.

$``(i)\Longrightarrow(iii)"$. Since
$$\begin{aligned}
 \mu_{\{b_{n}\},\{D_{n}\}}&=\delta_{b_{1}^{-1}\left\{0,1, \cdots,p_{1}-1\right\} t_{1}}\ast\delta_{b_{1}^{-1}p_{2}^{-1}\left\{0,1, \cdots,p_{2}-1\right\}t_{2} }
\ast\delta_{b_{1}^{-1}p_{2}^{-2}\left\{0,1, \cdots,p_{2}-1\right\}t_{2} }\ast\delta_{b_{1}^{-1}p_{2}^{-3}\left\{0,1, \cdots,p_{2}-1\right\} t_{2} }\ast \cdots
\\&=\delta_{b_{1}^{-1}\left\{0,1, \cdots,p_{1}-1\right\}t_{1}}\ast\frac{b_{1}}{t_{2}} \mathcal{L}_{[0,  \frac{t_{2}}{b_{1}}]}
\\&=\frac{b_{1}}{t_{2}p_{1}}\sum_{k=0}^{p_{1}-1}
\mathcal{L}_{[\frac{kt_{1}}{b_{1}},\frac{kt_{1}}{b_{1}}+\frac{t_{2}}{b_{1}}]}
 \\&=\frac{b_{1}}{t_{2}p_{1}}\left(\mathcal{L}_{ [0,  \frac{t_{2}}{b_{1}}]}+\mathcal{L}_{[\frac{t_{1}}{b_{1}},  \frac{t_{1}+t_{2}}{b_{1}} ]}+\mathcal{L}_{[ \frac{2t_{1}}{b_{1}},  \frac{2t_{1}+t_{2}}{b_{1}}]}+\cdots+\mathcal{L}_{[ \frac{(p_{1}-1)t_{1}}{b_{1}},  \frac{(p_{1}-1)t_{1}+t_{2}}{b_{1}}]}\right)
,\end{aligned}$$ where $\mathcal{L}_{[a, b]}$ is the Lebesgue measure restricting on the interval $[a, b]$.
According to $t_{2}\mid t_{1}$, we can let $t=\frac{t_{1}}{t_{2}}$. Then
$$ K:=K(b_{k}, D_{k})=\bigcup_{k=0}^{p_{1}-1}\left[kt\frac{t_{2}}{b_{1}}, (kt +1)\frac{t_{2}}{b_{1}} \right].$$
It is easy to check that $K\oplus J=\mathbb{R}$ and $ \mathcal{L}((K+j_{1})\cap(K+j_{2}))=0 $ for all $j_{1}\neq j_{2}\in J$, where $J= \frac{t_{2}}{b_{1}}\left(\left\{0,1,\cdots,t-1\right\}+tp_{1}\mathbb{Z}\right)$. Therefore, $K(b_{k},D_{k})$ is a  translation tile.

$``(iii)\Longrightarrow(i)"$. Since $K$
is a translation tile, then $t_{1}\geq t_{2}$ and there exists $J\subset\mathbb{R}$ such that $K\oplus J=\mathbb{R}$ and $ \mathcal{L}((K+j_{1})\bigcap(K+j_{2}))=0 $ for all $j_{1}\neq j_{2}\in J$.  Let $K_{i}=[\frac{it_{1}}{b_{1}},  \frac{it_{1}+t_{2}}{b_{1}}]$ for $i\in\{0,1,\cdots,p_{1}-1\}$, then $K=\bigcup_{i=0}^{p_{1}-1}K_{i}$. If $t_{1}=t_{2}$, we have $t_{2}\mid t_{1}$. If $t_{1}>t_{2}$, then there exists $r\in\mathbb{N} $ such that
 $\bigcup_{k=1}^{r}\left(K_{i_{k}}+j_{k}\right)=[\frac{t_{2} }{b_{1}}, \frac{t_{1}}{b_{1}}]$ for some $j_{k}\in J $, $i_{k} \in\{0,1,\cdots,p_{1}-1\}$ and $ \mathcal{L} \left(( K_{i_{k_{1}}}+j_{k_{1}})\bigcap(K_{i_{k_{2}}}+j_{k_{2}})\right)=0$ for any $k_{1}\neq k_{2}\in \{1,2,\cdots,r\}$. Hence, we have $$ \mathcal{L} \left(\bigcup_{k=1}^{r}\left(K_{i_{k}}+j_{k}\right) \right)=\sum_{k=1}^{r} \mathcal{L} \left(K_{i_{k}}+j_{k}\right)= \mathcal{L} \left(\left[\frac{t_{2} }{b_{1}}, \frac{t_{1}}{b_{1}}\right]\right).$$
This imply that $\frac{rt_{2}}{b_{1}}=\frac{t_{1}-t_{2}}{b_{1}}$, i.e., $t_{1}=(r+1)t_{2}$. Therefore, we have $t_{2}\mid t_{1}$.

In conclusion, we complete  the proof of Theorem \ref{thm(1.12)}.
\end{proof}
 The following  corollary  can be derived directly from Propositions \ref{prop(2.4)},  \ref{prop(3.3.2)} and Theorem  \ref{thm(1.12)}.
\begin{cor}\label{cor(1.1)}
 Given an integers  $N\geq 2$ and   integers $ |b_{k}|, p_{k}  \geq2, |t_{k}|\geq1$ for $1 \leq k\leq N$. Let  $|b_{k}|=p_{N}$ for all $k\geq N$ and $\mu_{\{b_k\},\{D_k\}}$ be defined by \eqref{eq(2.8)}, where
 $$D_{k}= \begin{cases}\left\{0,1, \cdots, p_{k}-1\right\}t_{k}, & if~1 \leq k\leq N-1 ;\\ \left\{0,1, \cdots, p_{N}-1\right\}t_{N} , & if~k\geq N.\end{cases}$$ Suppose that $\mu_{\{b_{k}\},\{D_k\}}$ is a spectral measure, then $t_{N}\mid t_{N-1}$.
\end{cor}
\begin{rem}
In  Corollary  \ref{cor(1.1)}, we can  get $t_{N}\mid t_{N-1}$, but we  can't  push  forward (i.e., we  can't get  $t_{N}\mid t_{N-2}$). The following example demonstrates this point well.
\end{rem}
\begin{exam}\label{ex(3.1)}
Given integers $ p_{1},  p_{2},p_{3}, t \geq2$, $t' \geq 1$ and an integer sequence $\{b_k\}_{k=1}^{\infty}$ with $b_{k} \geq2$.
Let$$
D_{k}= \begin{cases}\{0,1,\cdots,p_{1}-1\} , & if ~k=1;  \\  \{0,1,\cdots,p_{2}-1\}tt'   , & if~k= 2; \\ \{0,1,\cdots, p_{3}-1\}t , & if ~ k\geq 3, \end{cases}
$$where $\gcd(p_{2}, t')=1$, $ b_{2}=p_{2}t$ and $b_{k}=p_{3}$ for all $k \geq3$. Then $\mu_{\{b_{k}\},\{D_k\}}$ is a spectral measure.
\end{exam}
\begin{proof}
According to  Remark \ref{rem(2.3)}, the spectrality of $\mu_{\{b_{k}\},\{D_k\}}$ is independent of  $b_{1}$. Hence, we can assume  $b_{1}= p_{1}$. Let $\Lambda= \{0,1,\cdots,p_{1}-1\}+b_{1}\mathbb{Z}$.  Applying that $m_{D_{1}}(x)$ is a $ \mathbb{Z} $-periodic function and $m_{D_{2}}(x)$ is a $\frac{\mathbb{Z}}{t t'}$-periodic function, we have
$$ \begin{aligned}
  Q_{\mu_{\{b_{k}\},\{D_k\}},\Lambda} (x)&=\sum_{i=0}^{p_{1}-1}| m_{D_1} (\frac{i+x}{b_{1}} ) |^{2}  \sum_{\alpha\in \mathbb{Z}}| m_{D_2} (\frac{\alpha}{ b_{2}}+\frac{i +x}{b_{1}b_{2}} ) |^{2} |\widehat{\nu}_{>2} ( \frac{\alpha}{ b_{2}}+\frac{i +x}{b_{1}b_{2}} ) |^{2}
   \\&=\sum_{i=0}^{p_{1}-1}| m_{D_1} (\frac{i+x}{b_{1}} ) |^{2} \sum_{j=0}^{p_{2}-1}\sum_{\kappa\in \mathbb{Z}}| m_{D_2} (\frac{ j +\kappa p_{2} }{ b_{2}}+\frac{i +x}{b_{1}b_{2}} ) |^{2} |\widehat{\nu}_{>2} ( \frac{ j +\kappa p_{2}}{ b_{2}}+\frac{i +x}{b_{1}b_{2}} ) |^{2}
  \\& =\sum_{i=0}^{p_{1}-1}| m_{D_1} (\frac{i+x}{b_{1}} ) |^{2}\sum_{j=0}^{p_{2}-1} | m_{D_2} (\frac{ j }{ b_{2} }+\frac{i +x}{b_{1}b_{2}} ) |^{2} \sum_{\kappa\in \mathbb{Z}} |\widehat{\nu}_{>2} ( \frac{  \kappa }{ t}+ \frac{i +jb_{1}+x}{b_{1}b_{2}} ) |^{2}\end{aligned}$$
for all $x \in \mathbb{R}$.

 It is well known that $ \frac{\mathbb{Z}}{t }$  is a spectrum of  $\nu_{>2}$. From Proposition  \ref{prop(2.3)} $(ii)$,  we can know that $\sum_{\kappa\in \mathbb{Z}} |\widehat{\nu}_{>2} ( \frac{  \kappa }{ t}+ \frac{i +jb_{1}+x}{b_{1}b_{2}} ) |^{2}=1$.  Let $L_{1}=\{0,1,\cdots,p_{1}-1\}$ and $L_{2}=\{0,1,\cdots,p_{2}-1\}$.
According to $\gcd(p_{2}, t' )=1$, we have $t'\{0,1,\cdots,p_{2}-1\}\equiv \{0,1,\cdots,p_{2}-1\}({\rm mod}\ p)$. This means that $\frac{ (L_{2}-L_{2})\backslash\{0\}}{b_{2}}= \frac{ (t'L_{2}-t'L_{2})\backslash\{0\}}{p_{2}tt'}\subset  \frac{ \mathbb{Z}\setminus p \mathbb{Z}}{p_{2}tt'}$. Hence $ (b_{2}^{-1}D_2, L_{2})$ is  a compatible pair.  It is easy to verify that $ (b_{1}^{-1}D_1, L_{1})$ is also a compatible pair.   It follows from Proposition \ref{prop(2.2)} that $\sum_{j=0}^{p_{2}-1} | m_{D_2} (\frac{ j t'}{ p_{2}tt'}+\frac{i  +x}{b_{1}b_{2}} ) |^{2}=1$ and $\sum_{i=0}^{p_{1}-1}| m_{D_1} (\frac{i+x}{p_{1}} ) |^{2}=1$.
This shows that  $Q_{\mu_{\{b_{k}\},\{D_k\}},\Lambda} (x)=1$. Therefore, $\mu_{\{b_{k}\},\{D_k\}}$ is a spectral measure by  Proposition \ref{prop(2.3)} $(ii)$ .
 \end{proof}

\section{\bf   Proof of Theorem \ref{thm(1.13)} \label{sect.4}}
In this section, we focus on proving Theorem \ref{thm(1.13)}. Recall that, in Theorem \ref{thm(1.13)}, the finite pairs $ \{(b_{k},D_{k})\}_{k=1}^{m} :=\{(b_{k}, \{0,1,\cdots, p_{k}-1\}t_{k}) \}_{k=1}^{m}$, where $p_{k},t_{1},t_{2}, \cdots, t_{m}$ are pairwise coprime  integers for all $k\in\{1,2, \cdots, m\}$.
By Proposition \ref{prop(2.4)}, we can assume that all  $b_{k}$ and $t_{k}$ are positive integers. Let
$$m_{0}=\#\{k : t_{k}=1, 1\leq k\leq m\}.$$
When $m_{0}=m$,  Theorem \ref{thm(1.13)} can be obtained directly from Theorem \ref{thm(1)}. Hence,  we only need to consider the case that $m_{0}<m$.
From the definition of $\mu_\sigma$($\sigma \in\Sigma^{(m)}$), we  can assume that $ m_{0}\geq1$ (if $m_{0}=0$, we can add $D'=\{0,1,\cdots,p'-1\}$ into $\{D_k\}_{k=1}^m$ with  $p',t_{1},t_{2}, \cdots, t_{m}$ being pairwise coprime). Moreover, we also can assume that $t_{1}=t_{2}=\cdots=t_{m_{0}}=1$(by rearrange $\{D_k\}_{k=1}^m$ ). This means that $D_k$  can be written as
$$
D_{k}= \begin{cases}\left\{0,1, \cdots, p_{k}-1\right\}, & if~1\leq k\leq m_{0};\\ \left\{0,1, \cdots, p_{k}-1\right\}t_{k} , & if~k\in\{m_{0}+1, m_{0}+2,\cdots,m\}, \end{cases}
$$
where $ p_{k},t_{m_{0}+1}, t_{m_{0}+2}, \cdots, t_{m}$ are pairwise coprime for all $k\in\{1,2, \cdots, m\}$.

Let $ \Sigma^{(m)}=\{1,2, \cdots, m\}^{\mathbb{N}}$ be a symbolic space, we define  the metric
$$d(\sigma, \varsigma) =2^{- min\{n\geq1: \sigma_{n} \neq\varsigma_{n}\}}$$
for $\sigma = (\sigma_n)_{n=1}^{ \infty}$ and $\varsigma= (\varsigma_n)_{n=1}^{ \infty}\in \Sigma^{(m)}$, then $\Sigma^{(m)}$ is a compact metric space.
For a sequence $\{\sigma(i)\}_{i=1}^{ \infty}$ in $\Sigma^{(m)}$, we have that $ \sigma(i) $ converges to $\varsigma$ in $ \Sigma^{(m)}$ if and only if for each $n \geq 1$, there exists $i_{0} \geq1$ such that
$$\sigma_{1}(i)\sigma_{2}(i) \cdots\sigma_{n}(i) = \varsigma_{1}\varsigma_{2}\cdots\varsigma_{n}$$
for all $i\geq i_0$.
 For $\sigma\in\Sigma^{(m)}$ and the  finite  sequence pairs $\{(b _{k}, D_{k})\}_{k=1}^{m}$, we define the infinite convolution
\begin{equation} \label{eq(3.24)}
\mu_{\sigma}=\delta_{b_{\sigma_{1}}^{-1} D_{\sigma_{1}}} * \delta_{\left(b_{\sigma_{1}} b_{\sigma_{2}}\right)^{-1} D_{\sigma_{2}}} * \delta_{\left(b_{\sigma_{1}} b_{\sigma_{2}} b_{\sigma_{3}}\right)^{-1}D_{\sigma_{3}}}  * \cdots.
\end{equation}
Let $n \in \mathbb{N}$,  we write
$$
\mu_{\sigma, n}=\delta_{b_{\sigma_{1}}^{-1} D_{\sigma_{1}}} \ast  \delta_{\left(b_{\sigma_{1}} b_{\sigma_{2}}\right)^{-1} D_{\sigma_{2}}}\ast \cdots * \delta_{\left(b_{\sigma_{1}} b_{\sigma_{2}} \cdots b_{\sigma_{n}}\right)^{-1}D_{\sigma_{n}}}.
$$
and
\begin{equation} \label{eq(3.2444)}
\nu_{\sigma, >n}=\delta_{b_{\sigma_{n+1}}^{-1} D_{\sigma_{n+1}}} \ast \delta_{\left(b_{\sigma_{n+1}} b_{\sigma_{n+2}}\right)^{-1} D_{\sigma_{n+2}}}\ast \delta_{\left(b_{\sigma_{n+1}} b_{\sigma_{n+2}} b_{\sigma_{n+3}}\right)^{-1}D_{\sigma_{n+3}}} \ast\cdots .
\end{equation}
For any  fixed $\sigma\in \Sigma^{(m)}$, it is clear that $\mu_{\sigma, n}$ converges weakly to  $\mu_{\sigma}$ by \eqref{eq(2.82)}.

Next, we will give the definition of integral periodic zero set of measure $\mu$.
\begin{defi}\label{de(4.1)}  Let $\mu$ be a Borel probability measure in $\mathbb{R}$.  Its integral periodic zero set is defined by
$$
 \mathrm{Z}(\mu)=\{\xi \in \mathbb{R}: \widehat{\mu}(\xi+k)=0~~\text{for all}~ k \in \mathbb{Z}\} .
$$
\end{defi}
In \cite {LMW-2022-2},  Li et al.  given the following  theorem, which plays a key role in the proof of  Theorem \ref{thm(1.13)}.
\begin{thm}\label{thm(4.1)}\cite[Theorem 1.1]{LMW-2022-2}
 Given a sequence of admissible pairs $\left\{\left(b_{k}, D_{k}\right)\right\}_{k=1}^{\infty}$ in $\mathbb{R}$. Suppose that the infinite convolution $\mu$ defined in \eqref{eq(2.8)} exists. Let the sequence $\{\nu_{>k}\}_{k=1}^{\infty}$ be defined in \eqref{eq(2.90)}. If there exists a subsequence $\{\nu_{>k_{j}}\}_{j=1}^{\infty}$ which converges weakly to $\nu$, and $ \mathrm{Z} (\nu)=\emptyset$, then $\mu$ is a spectral measure.
 \end{thm}
In general, it is hard to  prove that  the integral periodic zero set  is an empty set. The following three propositions give some sufficient conditions for integral periodic zero set to be an empty set.
\begin{prop}\label{thm(4.2)}\cite{LMW-2022-2}
 Let the sequence of admissible pairs $\{(b_{k}, D_{k})\}_{k=1}^{\infty}$ be chosen from a finite set of admissible pairs in $\mathbb{R}$, and let $\mu $ be defined by \eqref{eq(2.8)}. Suppose that
$$
\operatorname{gcd}\left(\bigcup_{k=n}^{\infty}(D_{k}-D_{k})\right)=1
$$
for each $n \geq 1$, then $\mathrm{Z}(\mu)=\emptyset$.
\end{prop}
Define
\begin{equation} \label{eq(3.24000)}A=\{1,2,\cdots,m_{0}\} \end{equation}
and
\begin{equation} \label{eq(3.24001)}  B=\{m_{0}+1,m_{0}+2,\cdots,m\}.\end{equation}
If $\sigma=i^{\infty}$ with $i\in B$ and $p_{i}=b_{i}$, then $ \frac{1}{t_{i}}+\mathbb{Z}\subset \mathcal{Z}(\widehat{\mu}_{\sigma})=\bigcup_{k=0}^{\infty} \frac{p_{i}^{k}(\mathbb{Z} \backslash p_{i} \mathbb{Z})}{t_{i}}=\frac{\mathbb{Z}}{t_{i}}$. This means that $\mathrm{Z}(\mu_{\sigma})\neq\emptyset$. However, if $p_{i}\mid b_{i} $ and $p_{i}\neq b_{i}$,
we can combine Theorem  4.1 and Lemma 6.3  in \cite{LMW-2022-2} to get that $\mathrm{Z}(\mu_{\sigma})=\emptyset$.
\begin{prop}\label{thm(4.3)}
Let $\sigma=i^{\infty}$ with $i\in B$, and let $B$, $\mu_{\sigma}$ be defined by  \eqref{eq(3.24001)} and \eqref{eq(3.24)}, respectively. Suppose  that $p_{i}\mid b_{i} $ and $p_{i}\neq b_{i} $,  then $\mathrm{Z}(\mu_{\sigma})=\emptyset$.
\end{prop}

For some special $\sigma\in \Sigma^{(m)}$, according to the characteristics of the zeros of  Fourier transform of $\mu_{\sigma}$, we can prove that the integral periodic zero set of $\mu_{\sigma}$  is an empty set.
\begin{prop}\label{lem(4.1)}Let $\sigma= (\sigma_n)_{n=1}^{ \infty}\in \Sigma^{(m)}$ and  $\mu_{\sigma}$ be defined by  \eqref{eq(3.24)}. Suppose that $p_{\sigma_n}\mid b_{\sigma_n}$ for  all $n\geq1$ and one of the following holds:
\begin{enumerate}[(i)]
  \item $\sigma=i\xi$ with $i\in A$ and $\xi\in\Sigma^{(m)}$;
  \item $\sigma=j\eta$ with $j\in B$ and $\eta\in\{1,2,\cdots,j-1,j+1,\cdots,m\}^{\mathbb{N}}$,
\end{enumerate}
where $A$ and $B$ are defined by \eqref{eq(3.24000)} and \eqref{eq(3.24001)}, respectively. Then $\mathrm{Z}(\mu_{\sigma})=\emptyset$.\end{prop}
\begin{proof}
We will prove the proposition by contradiction. Suppose that $\mathrm{Z} (\mu_{\sigma})\neq\emptyset$, then there exists $\xi_{0} \in \mathrm{Z}\left(\mu_{\sigma}\right)$ such that $\widehat{\mu}_{\sigma}\left(\xi_{0}+k\right)=0$ for all $k \in \mathbb{Z}$. According to $\widehat{\mu}_{\sigma} ( 0 )=1$, this is easy to verify that $\xi_{0}+ \mathbb{Z}\subset \mathcal{Z}(\widehat{\mu}_{\sigma}) \backslash\mathbb{Z}$. Since $p_{\sigma_n}\mid b_{\sigma_n}$ for  all $n\geq1$,
it follows from \eqref{eq(2.5)} that
\begin{equation}\label{eq(4.2500)}\xi_{0}+ \mathbb{Z}\subset \mathcal{Z}(\widehat{\mu}_{\sigma}) \backslash\mathbb{Z}\subset  \frac{b(\mathbb{Z} \backslash p_{\sigma_{1}} \mathbb{Z})}{t_{\sigma_{1}}}  \bigcup \left(\bigcup_{\alpha\in B}\frac{p_{\sigma_{1}}  \mathbb{Z}}{t_{\alpha}}\right) \subset \bigcup_{\alpha\in B} \frac{ \mathbb{Z}}{t_{\alpha}},
\end{equation} where $b\in\mathbb{N} $ and $b = \frac{b_{\sigma_{1}}}{ p_{\sigma_{1}}}$.
We first prove the following two claims.

\noindent\emph{{\bf Claim 1.}  Suppose that there exist $k_{0}\in\mathbb{Z}$ and $\alpha\in B$ such that $\xi_{0} +k_{0}\in\frac{\mathbb{Z}}{t_{\alpha}}$, then $\xi_{0} +\mathbb{Z}\subset\frac{\mathbb{Z}}{t_{\alpha}}$.}	
\begin{proof}[Proof of Claim 1]
 Suppose, on the contrary, that there exists $\tilde{k}_{0} \in\mathbb{Z}$  such that $\xi_{0} +\tilde{k}_{0}\notin\frac{\mathbb{Z}}{t_{\alpha}}$. Since $\xi_{0}+ \mathbb{Z} \subset \bigcup_{\alpha\in B} \frac{ \mathbb{Z}}{t_{\alpha}}$, there exists  $\tilde{\alpha}\in B\backslash\{\alpha\}$ such that  $\xi_{0}+\tilde{k}_{0}\in\frac{\mathbb{Z}}{t_{\tilde{\alpha}}}$.
Then
$$
\xi_{0}+k_{0 }=\frac{k_{1} }{t_{\alpha}}\quad\text{and}\quad  \xi_{0}+\tilde{k}_{0 }=\frac{\tilde{k}_{1}}{t_{\tilde{\alpha}}}
$$
for some $k_{1}, \tilde{k}_{1} \in\mathbb{Z}$. The above equations imply that $\frac{t_{\alpha} \tilde{k}_{1}}{t_{\tilde{\alpha}}} = k_{1}+t_{\alpha}(\tilde{k}_{0 }- k_{0 })\in \mathbb{Z}$.
This contradicts with $\gcd(t_{a}, t_{\tilde{a}})=1$ and $\frac{\tilde{k}_{1}}{t_{\tilde{a}}}  \notin \mathbb{Z}$. The claim follows.
\end{proof}

 \noindent\emph{{\bf Claim 2.} $(\xi_{0}+ \mathbb{Z})\bigcap (\bigcup_{\alpha \in B\backslash\{\sigma_{1}\}} \frac{p_{\sigma_{1}} \mathbb{Z}}{t_{\alpha}})=\emptyset$.}	
\begin{proof}[Proof of Claim 2]  Suppose, on the contrary, that there exists
$\alpha_{0}\in B\backslash\{\sigma_{1}\}$ such that $(\xi_{0}+ \mathbb{Z})\bigcap \frac{p_{\sigma_{1}} \mathbb{Z}}{t_{\alpha_{0}}}\neq\emptyset$.
It follows from Claim 1 and \eqref{eq(4.2500)} that $\xi_{0}+ \mathbb{Z}\subset \frac{p_{\sigma_{1}} \mathbb{Z}}{t_{\alpha_{0}}}$. Hence $\xi_{0} =\frac{p_{\sigma_{1}}k_{2} }{t_{\alpha_{0}}}$ and $\xi_{0}+ 1=\frac{p_{\sigma_{1}}\tilde{k}_{2} }{t_{\alpha_{0}}}$  for some $k_{2}, \tilde{k}_{2} \in \mathbb{Z}$.
 By a simple calculation, we obtain that $t_{\alpha_{0}}=p_{\sigma_{1}}(\tilde{k}_{2}-k_{2} )$, which contradicts with $\gcd(p_{\sigma_{1}}, t_{\alpha_{0}})=1$ and $p_{\sigma_{1}}\geq2$.  This completes the proof of Claim 2.
\end{proof}
$(i)$. Since $\sigma=i\xi$ with $i\in A$ and $\xi\in\Sigma^{(m)}$, we infer from \eqref{eq(2.5)} that
$\xi_{0}+ \mathbb{Z} \subset \mathcal{Z}(\widehat{\mu}_{\sigma}) \backslash\mathbb{Z}\subset\bigcup_{\alpha\in B} \frac{p_{i}  \mathbb{Z}}{t_{\alpha}}.$ This contradicts with Claim 2. Therefore, we have $\mathrm{Z} (\mu_{\sigma})=\emptyset$.

 $(ii)$. Since $\sigma=j\eta$ with $j\in B$ and $\eta\in\{1,2,\cdots,j-1,j+1,\cdots,m\}^{\mathbb{N}}$, it follows from \eqref{eq(2.5)} that
\begin{equation} \label{eq(4.24)}
\xi_{0}+ \mathbb{Z} \subset \mathcal{Z}(\widehat{\mu}_{\sigma}) \backslash\mathbb{Z} \subset \frac{b(\mathbb{Z} \backslash p_{j} \mathbb{Z})}{t_{j}}  \bigcup \left(\bigcup_{\alpha\in B\backslash\{j\}}\frac{p_{j}  \mathbb{Z}}{t_{\alpha}}\right).
\end{equation}
 Note that $\gcd(p_{j}, t_j)=1$, then there exist $\iota_{1},\iota_{2}\in \mathbb{Z}$ such that $\iota_{1} t_{j}+\iota_{2}p_{j}=1$.  Applying \eqref{eq(4.24)} and Claim 2,  we conclude $ \xi_{0}+ \mathbb{Z}\subset\frac{ b (\mathbb{Z} \backslash p_{j} \mathbb{Z} )}{t_{j}}$. Hence, there exists  $k_{3} \in \mathbb{Z} \backslash p_{j} \mathbb{Z}$ such that $\xi_{0} =\frac{ b k_{3} }{t_j}$ and $\xi_{0}-\iota_{1} bk_{3}=\frac{b \tilde{k}_{3} }{t_j}$ for some $\tilde{k}_{3} \in \mathbb{Z} \backslash p_{j} \mathbb{Z}$.
Then $\frac{ k_{3} }{t_j} -\iota_{1} k_{3}=\frac{ \tilde{k}_{3} }{t_j}$, i.e., $ k_{3}-\iota_{1} k_{3}\frac{1-\iota_{2}p_{j}}{\iota_{1}} =  \tilde{k}_{3} $. By a simple calculation, we obtain $\tilde{k}_{3} \in p_{j} \mathbb{Z}$, which contradicts with $\tilde{k}_{3} \in \mathbb{Z} \backslash p_{j} \mathbb{Z}$. This means that $\mathrm{Z} (\mu_{\sigma})=\emptyset$.

The proof of  Proposition \ref{lem(4.1)}  is now completed.
\end{proof}

\begin{lem}\label{lem(4.5)}
Given a sequence $\{\sigma(i)\}_{i=1}^{ \infty}\subset\Sigma^{(m)}$ and $\sigma \in \Sigma^{(m)}$. Suppose that  $\sigma(i)$ converges to $\sigma$, then  $\mu_{\sigma(i)}$ converges weakly  to $\mu_{\sigma  }$.
\end{lem}
\begin{proof}
 For any $f \in C_{b}(\mathbb{R})$ and $\varepsilon > 0$,  it follows from  \cite[Lemma 5.1] {LMW-2022-2}  that  there exists $q_{0}\geq1$ such that
\begin{equation} \label{eq(4.444)}\left|\int_{\mathbb{R}} f(x) d\mu_{\sigma, q_{0}}(x)-\int_{\mathbb{R}} f(x) d\mu_{\sigma}(x)\right|<\frac{\varepsilon}{2}
\end{equation}
and
\begin{equation} \label{eq(4.445)}
\left|\int_{\mathbb{R}} f(x) d\mu_{\sigma(i), q_{0}}(x)-\int_{\mathbb{R}} f(x) d \mu_{\sigma(i)}(x)\right|<\frac{\varepsilon}{2}
\end{equation} for $i\geq1$.
Since  ${\sigma(i)}$ converges to $\xi$, there exists $i_{0}\geq1$ such that
$\mu_{\sigma(i), q_{0}}=\mu_{\sigma, q_{0}}$ for $i\geq i_{0}$.
We conclude from \eqref{eq(4.444)} and \eqref{eq(4.445)} that
$$ \left|\int_{\mathbb{R}} f(x) d\mu_{\sigma(i)}(x)-\int_{\mathbb{R}} f(x) d\mu_{\sigma}(x)\right|<\varepsilon$$
for $i\geq i_{0}$. Therefore, $ \mu_{\sigma(i)} $ converges weakly to $\mu_{\sigma}$.
 \end{proof}
For a positive integer $l$, we define
\begin{equation} \label{eq(4.244)}\Xi_{l}=\left\{i_{1}i_{2}\cdots i_{l}j^{\infty}\in\Sigma^{(m)}:  i_{l}\neq j,  j\in B\right\},
\end{equation}
\begin{equation} \label{eq(4.245)}\Xi_{l}^{(1)}=\left\{i_{1}i_{2}\cdots i_{l}j^{\infty}\in\Sigma^{(m)}:  i_{l}\neq j,  j\in B, b_{j}= p_{j}\right\}
\end{equation}
and
\begin{equation} \label{eq(4.246)}\Xi_{l}^{(2)}=\left\{i_{1}i_{2}\cdots i_{l}j^{\infty}\in\Sigma^{(m)}:  i_{l}\neq j,  j\in B, b_{j}\neq p_{j}\right\},
\end{equation}
where $B$ is defined by  \eqref{eq(3.24001)}. It is obvious that $\Xi_{l}=\Xi_{l}^{(1)}\bigcup\Xi_{l}^{(2)}$.
\begin{prop}\label{prop(4.1)}
For $\sigma\in\Sigma^{(m)}$, let $\mu_{\sigma}$  and $ \Xi_{l} $  be defined by  \eqref{eq(3.24)} and  \eqref{eq(4.244)}  respectively.  Suppose that $\sigma\notin\bigcup_{l=1}^{\infty}\Xi_{l} $ and $p_{\sigma_n}\mid b_{\sigma_n}$ for  all $n\geq1$, then the infinite convolution $\mu_{\sigma}$ is a spectral measure.
\end{prop}
 \begin{proof}
Since $\gcd(p_{\sigma_n}, t_{\sigma_n})=1$ and $p_{\sigma_n}\mid b_{\sigma_n}$ for all $n\geq1$, it follows from Proposition \ref{prop(2.5)} that $\left(b_{\sigma_n}, D_{\sigma_n}\right)$  are admissible pairs for all $n\geq1$. We define
$$
\varrho^{n}(\sigma)=\sigma_{n+1}\sigma_{n+2} \sigma_{n+3}\cdots
$$
for $n\geq1$, and let $A$ and $B$ are defined by \eqref{eq(3.24000)} and \eqref{eq(3.24001)}, respectively. It is easy to verify that $\nu_{\sigma,>n}=\mu_{\varrho^{n}(\sigma)}$ for $n\geq1$, where $\nu_{\sigma,>n}$ is defined by \eqref{eq(3.24)}. We divide $\sigma$ into the following three  cases to discuss.

Case 1: $\sigma=i_{0}^{\infty}$ for some $i_{0}\in B$.

According to Lemma  \ref{lem(2.3)} and  Theorem \ref{thm(1)}, we can get that $\mu_{\sigma}$ is a spectral measure.

Case 2:  There exists $i_{0}\in A$ such that the symbol $`` i_{0}"$  occurs infinitely many times in $\sigma$.

 Let
$$
\{k_{1}, k_{2}, k_{3},\cdots\}=\{n \geq 1: \sigma_{n}=i_{0}\}
$$ and $k_{n}<k_{n+1}$ for all $n\geq1$. It is clear that the sequence $\{{\varrho^{k_{n}-1} (\sigma)}\}_{n=1}^{\infty} \subset\Sigma^{(m)}$  and $ \varrho^{k_{n}-1} (\sigma)=i_{0}\eta_{n}$ for some $\eta_{n}\in \Sigma^{(m)}$. By the compactness of $\Sigma^{(m)}$,  there exist  a subsequence $\{n_j\}_{j=1}^{\infty}$ in $\mathbb{N}$ and  $\tilde{\eta } \in\Sigma^{(m)}$ such that $ \varrho^{k_{n_{j}}-1} (\sigma) $ converges to $\eta:=i_{0}\tilde{\eta}$. Since $\{\nu_{\sigma,>(k_{n_{j}}-1)}\}_{j=1}^{\infty}=\{\mu_{\varrho^{k_{n_{j}}-1}(\sigma)}\}_{j=1}^{\infty}$, it follows from  Lemma \ref{lem(4.5)} and Proposition \ref{lem(4.1)} $(i)$ that $\nu_{\sigma,>(k_{n_{j}}-1)}$ converges weakly to $\mu_{\eta}$ and $\mathrm{Z}(\mu_{\eta})=\emptyset$.  In view of Theorem \ref{thm(4.1)}, $\mu_{\sigma}$ is a spectral measure.

Case 3: There exist  $i, j \in B $  and $i\neq j$  such that   two symbols $``i "$ and $``j  "$  occur infinitely many times in $\sigma$.

For $\kappa\in\{i, j\}$, let
$$
\{h_{1}^{(\kappa)}, h_{2}^{(\kappa)}, h_{3}^{(\kappa)}, \cdots\}=\{n \geq 1: \sigma_{n}=\kappa\}
$$
and $h_{n}^{(\kappa)}<h_{n+1}^{(\kappa)}$ for all $n\geq1$.

Suppose that $\varlimsup\limits _{n \rightarrow \infty}\left(h_{n+1}^{(\kappa)}-h_{n}^{(\kappa)}\right)<\infty$  for all $\kappa\in\{i, j\}$.
Since the sequence $\left\{\varrho^{n} (\sigma)\right\}_{n=1}^{\infty}\subset \Sigma^{(m)} $ and $\Sigma^{(m)}$ is a compact space, there exist a subsequence $\{n_j\}_{j=1}^{\infty}$ in $\mathbb{N}$ and $\xi\in\Sigma^{(m)}$ such that $ \varrho^{n_{j}} (\sigma) $ converges
to $\xi$. It is easy to know that  symbols  $``i "$ and $``j  "$  occur  infinitely many times in $\xi$. Since $\operatorname{gcd}\left(\bigcup_{n=k}^{\infty}(D_{\xi_{n}}-D_{ \xi_n})\right)=1$ for each $k \geq 1$, it follows from  Lemma \ref{lem(4.5)} and Proposition  \ref{thm(4.2)} that $\nu_{\sigma,>{n_{j}}}$ converges weakly to $\mu_{\xi}$ and $\mathrm{Z}(\mu_{\xi})=\emptyset$.
 Combining this with Theorem \ref{thm(4.1)}, we obtain that $\mu_{\sigma}$ is a spectral measure.

Suppose that $\varlimsup\limits_{n \rightarrow \infty}\left(h_{n+1}^{(\kappa)}-h_{n}^{(\kappa)}\right)=\infty$ for some $\kappa\in\{i, j\}$.
We may find a subsequence $\{n_{j}\}_{j=1}^{\infty}$ in $\mathbb{N}$ such that $\lim\limits_{j \rightarrow \infty} \left(h_{n_{j}+1}^{(\kappa)}-h_{n_{j}}^{(\kappa)}\right )=\infty.$  It is clear that $\varrho^{h_{n_{j}}^{(\kappa)}-1} (\sigma)=\kappa\zeta_{j}$ for some $\zeta_{j}\in \Sigma^{(m)}$. Hence, there exist a subsequence $\{\tilde{n}_j\}_{j=1}^{\infty}$ in $\{n_{j}\}_{j=1}^{\infty}$ and $\vartheta:=\kappa \zeta $ with $ \zeta \in \Sigma^{(m)}$ such that ${\varrho^{h_{\tilde{n}_{j}}^{(\kappa)}-1} (\sigma)}$ converges to $\vartheta$.  Since $\lim\limits_{j \rightarrow \infty} \left(h_{n_{j}+1}^{(\kappa)}-h_{n_{j}}^{(\kappa)}\right )=\infty$, we have $\zeta \in \{1,2,\cdots,\kappa-1,\kappa+1, \cdots, m\}^{\mathbb{N}}$. It follows from  Lemma \ref{lem(4.5)}, Proposition \ref{lem(4.1)} $(ii)$ and Theorem \ref{thm(4.1)} that $\mu_{\sigma}$ is a spectral measure.

 The proof of Proposition \ref{prop(4.1)} is completed.
\end{proof}
 \begin{prop}\label{prop(4.6)}
For $\sigma\in\Sigma^{(m)}$, let $\mu_{\sigma}$, $\Xi_{l}$  and $ \Xi_{l}^{(2)}$ be defined by  \eqref{eq(3.24)}, \eqref{eq(4.244)}  and \eqref{eq(4.246)}, respectively.  Suppose that  $\sigma\in\bigcup_{l=1}^{\infty}\Xi_{l}$ and $p_{\sigma_n}\mid b_{\sigma_n}$ for  all $n\geq1$, then $\mu_{ \sigma}$ is a spectral measure if and only if $\sigma\in\bigcup_{l=1}^{\infty}\Xi_{l}^{(2)}$.
\end{prop}
 \begin{proof}
We first prove the necessity. Suppose,  on the contrary, that  $\sigma\notin\bigcup_{l=1}^{\infty}\Xi_{l}^{(2)}$, then there exists $l_{1}\geq1$ such that $\sigma\in \Xi_{l_{1}}^{(1)}=\left\{i_{1}i_{2}\cdots i_{l_{1}}j^{\infty}\in\Sigma^{(m)}:  i_{l_{1}}\neq j,  j\in B, b_{j}= p_{j}\right\}$. Hence, $\sigma=\sigma_{1}\sigma_{2}\cdots\sigma_{l_{1}}j_{0}^{\infty}$ and $ b_{j_{0}}= p_{j_{0}}$  for some  $j_{0}\in B\backslash\{ \sigma_{l_{1}}\}$. Since  $\gcd( t_{j_{0}} , t_{\sigma_{l_{1}}})=1$ and $t_{j_{0}}\neq1 $, it follows from  Theorem \ref{thm(1.12)}
that $$ \nu_{\sigma,>(l_{1}-1)}=\delta_{b_{\sigma_{l_{1}}}^{-1}D_{\sigma_{l_{1}}}}\ast \delta_{b_{\sigma_{l_{1}}}^{-1}p_{j_{0}}^{-1}D_{j_{0}}}\ast\delta_{b_{\sigma_{l_{1}}}^{-1}p_{j_{0}}^{-2}D_{j_{0}}}  \ast\delta_{b_{\sigma_{l_{1}}}^{-1}p_{j_{0}}^{-3}D_{j_{0}}} \ast\cdots$$
 is not a spectral measure. According to Propositions \ref{prop(3.3.2)}, we have $\mu_{ \sigma}$ is not a spectral measure, which is a contradiction. Hence the necessity follows.

Next, we prove the  sufficiency.  Since $\sigma\in\bigcup_{l=1}^{\infty}\Xi_{l}^{(2)}$, there exist $l_{2}\geq1$ and $j_{0}\in B$ such that $\sigma =\sigma_{1}\sigma_{2}\cdots \sigma_{l_{2}}j_{0}^{\infty}$ and $b_{j_{0}}\neq p_{j_{0}}$ and
 $$ \nu_{\sigma,>k}=\delta_{b_{j_{0}}^{-1}D_{j_{0}}}\ast\delta_{b_{j_{0}}^{-2}D_{j_{0}}}  \ast\delta_{b_{j_{0}}^{-3}D_{j_{0}}} \ast\cdots $$  for all $k\geq l_{2}$. From Lemma \ref{lem(4.5)}, we can easily know that $ \nu_{\sigma,>k}$ converges weakly to $\tilde{\nu}:=\delta_{b_{j_{0}}^{-1}D_{j_{0}}}\ast\delta_{b_{j_{0}}^{-2}D_{j_{0}}}  \ast\delta_{b_{j_{0}}^{-3}D_{j_{0}}} \ast\cdots$. Since $  p_{j_{0}}\mid  b_{j_{0}}$ and $p_{j_{0}}\neq b_{j_{0}}$,  it follows from  Proposition \ref{thm(4.3)} that $\mathrm{Z}(\tilde{\nu})=\emptyset$. According to Proposition \ref{prop(2.5)}, we have  $\left(b_{\sigma_n}, D_{\sigma_n}\right)$  are admissible pairs for all $n\geq1$. Then $\mu_{\sigma}$ is  a spectral measure by Theorem \ref{thm(4.1)}. This proves the sufficiency, and hence the proof is completed.
\end{proof}
Based on the above preparations, now we can prove Theorem \ref{thm(1.13)}.
\begin{proof}[\bf Proof of  Theorem \ref{thm(1.13)}]
According to Proposition \ref{prop(2.4)}, we can assume that $b_{k}\geq2$ and $t_{k}\geq1$ for  $k\in\{1,2, \cdots, m\}$.  We also can assume that  $p_{\sigma_1}\mid b_{\sigma_1}$ by Remark \ref{rem(2.3)}. Therefore, the proof of Theorem \ref{thm(1.13)} can be directly obtained from Proposition \ref{prop(4.1)} and Proposition \ref{prop(4.6)} and Theorem \ref{thm(1.1)}.
\end{proof}
At the end of this paper, we give some examples which related to our main results.%
\begin{exam}\label{ex(4.4)}
Given integers $p_{1}, p_{2}\geq2$ and  integer sequences $\{b_{k}\}_{k=1}^{\infty}$, $\{t_{2k}\}_{k=1}^{\infty}$ with $b_{k} \geq2$ and $t_{2k} \geq1$. Let $\mu_{\{b_k\},\{D_k\}}$ be defined by \eqref{eq(2.8)}, where
$$
D_{k}= \begin{cases}\{0,1,\cdots, p_{1}-1\} , & if ~k\in 2\mathbb{Z}+1;  \\  \{0,1,\cdots,p_{2}-1\}t _{k}, & if~k\in 2\mathbb{Z}.\end{cases}
$$
Suppose that $b_{2k}=p_2t_{2k}$ and the sequence $\left\{t_{2k}\right\}_{k=1}^{\infty}$ is bounded. Then $\mu_{\{b_k\},\{D_k\}}$ is a spectral measure if and only if $p_{1}\mid b_{2k+1}t_{2k}$ for  $k\geq1$.
\end{exam}
\begin{proof}
The necessity follows directly from Proposition \ref{prop(3.4.2)}.
Now we are devoted to proving the sufficiency. By  Lemma  \ref{lem(2.3)}, we  know that  $\mu_{\{b_k\},\{D_k\}}$ is a spectral measure if and only if  $\mu_{\{b_k\},\{p_{2} D_{k}\}}$ is a spectral measure, where
$$\begin{aligned}
\mu_{\{b_k\},\{p_{2}D_k\}}&=(\delta_{b_{1}^{-1}p_{2}D_1}\ast\delta_{b_{1}^{-1}b_{2}^{-1}p_{2}D_2})
\ast(\delta_{b_{1}^{-1}b_{2}^{-1} b_{3}^{-1}p_{2}D_3}\ast\delta_{b_{1}^{-1}b_{2}^{-1} b_{3}^{-1} b_{4}^{-1}p_{2}D_4})\ast\cdots
\\&\cdots\ast(\delta_{b_{1}^{-1}b_{2}^{-1}\cdots  b_{2k}^{-1}b_{2k+1}^{-1}p_{2}D_{2k+1}}\ast\delta_{b_{1}^{-1}b_{2}^{-1}\cdots b_{2k}^{-1}b_{2k+1}^{-1}b_{2k+2}^{-1}p_{2}D_{2k+2}})\ast\cdots.
\end{aligned}$$
By a simple calculation, we obtain
 $$\begin{aligned}
 \delta_{b_{1}^{-1}b_{2}^{-1}\cdots b_{2k+1}^{-1}p_{2}D_{2k+1}}\ast\delta_{b_{1}^{-1}b_{2}^{-1}\cdots b_{2k+1}^{-1}b_{2k+2}^{-1}p_{2}D_{2k+2}}
 &=\delta_{b_{1}^{-1}b_{2}^{-1}\cdots b_{2k}^{-1} b_{2k+1}^{-1}\{0,1,\cdots, p_{1}-1\}p_{2}}\ast\delta_{b_{1}^{-1}b_{2}^{-1}\cdots b_{2k}^{-1} b_{2k+1}^{-1} \{0,1,\cdots, p_{2}-1\}}
 \\&= \delta_{b_{1}^{-1}b_{2}^{-1}\cdots b_{2k}^{-1} b_{2k+1}^{-1}\{0,1,\cdots, p_{1}p_{2}-1\}}
  \end{aligned}$$ for any $k\geq1$.
 Then
 $$\begin{aligned}
\mu_{\{b_k\},\{p_{2}D_k\}}&= \delta_{b_{1}^{-1}\{0,1,\cdots, p_{1}p_{2}-1\}}
\ast\delta_{b_{1}^{-1}(b_{2}^{-1} b_{3}^{-1})\{0,1,\cdots, p_{1}p_{2}-1\}}\ast\delta_{b_{1}^{-1}(b_{2}^{-1} b_{3}^{-1})(b_{4}^{-1} b_{5}^{-1})\{0,1,\cdots, p_{1}p_{2}-1\}}\ast\cdots\\& \cdots\ast
 \delta_{b_{1}^{-1}(b_{2}^{-1}b_{3}^{-1})\cdots  (b_{2k}^{-1}b_{2k+1}^{-1})\{0,1,\cdots, p_{1}p_{2}-1\}} \ast\cdots .
\end{aligned}$$
Since  $p_{1}\mid b_{2k+1}t_{2k}$ and $b_{2k}=p_2t_{2k}$ for $k\geq1$, we have $p_{1}p_2\mid b_{2k}b_{2k+1}$ for $k\geq1$. This together with  Theorem \ref{thm(1)} shows that $\mu_{\{b_k\},\{p_{2}D_k\}}$ is a spectral measure. Therefore, $\mu_{\{b_k\},\{D_k\}}$ is a spectral measure.
\end{proof}
\begin{exam}
Given integers $|b|, N\geq2 $ and  two coprime integers  $ p, |t|\geq2 $. Let $\mu_{b,\{D_k\}}$ be defined by \eqref{eq(2.8)}, where
$$
D_{k}= \begin{cases}\{0,1,\cdots, p-1\} , & if~k< N; \\  \{0,1,\cdots, p-1\}t , & if~k\geq N.  \end{cases}
$$
Then $\mu_{b,\{D_k\}}$ is a spectral measure if and only if  $ p\mid b$ and $| b |\neq p$.
\end{exam}
\begin{proof}
The proof follows directly from Theorem \ref{thm(1.13)}.
\end{proof}

\end{document}